\documentclass[10pt]{ijnam}
\usepackage{epsfig}
\usepackage{amsmath}
\usepackage{pst-plot}
\usepackage{amssymb}
\usepackage{amsfonts}
\usepackage{euscript}
\usepackage{eufrak}
\usepackage{oldgerm}

\renewcommand{\theequation}{\thesection.\arabic{equation}}
\newcommand{\se}{\setcounter{equation}{0}}

\newtheorem{theorem}{Theorem}[section]
\newtheorem{lemma}{Lemma}[section]
\newtheorem{remark}{Remark}[section]

\newcommand{\cA}{\mathcal{A}}
\newcommand{\cB}{\mathcal{B}}
\newcommand{\cT}{\mathcal{T}}
\newcommand{\cE}{\mathcal{E}}

\newcommand{\G}{G}
\newcommand{\bG}{\hat{G}}
\newcommand{\hh}{\hat{\theta}}
\newcommand{\rv}{V_h}

\newcommand{\vv}{U_h^*}

\newcommand{\pp}{{S}}
\newcommand{\bpp}{\hat{S}}

\newcommand{\beh}{{\bar\epsilon}_h}
\newcommand{\beA}{{\bar\epsilon}_A}
\newcommand{\beB}{{\bar\epsilon}_B}
\newcommand{\tAh}{\tilde{A}_h}
\newcommand{\tBh}{\tilde{B}_h}

\newcommand{\ph}{{\partial}_t}
\newcommand{\pc}{{\delta}_t}

\newcommand{\pb}{\bar{\partial}_t}
\newcommand{\ps}{{\partial}_{t}^2}
\newcommand{\dt}{k}

\copyrightinfo{2004}{}

\begin{document}

\title [FVEM for Hyperbolic PIDE]
{A priori Error Estimates for Finite Volume Element Approximations to Second Order Linear Hyperbolic 
Integro-Differential Equations}
\author[S. Karaa]{Samir Karaa}
\address{
Department of Mathematics and Statistics, Sultan Qaboos University, 
P. O. Box 36, Al-Khod 123, Muscat, Oman
}
\email{skaraa@squ.edu.om}

\author[A. K. Pani]{Amiya K. Pani}
\address{
Department of Mathematics, Industrial Mathematics Group, 
Indian Institute of Technology,
Bombay, Powai, Mumbai-400076, India
}
\email{akp@math.iitb.ac.in}

%




\abstract{In this paper, both semidiscrete and completely discrete finite volume element methods  
(FVEMs) are analyzed for approximating solutions of  a class of linear hyperbolic integro-differential equations in a two-dimensional convex polygonal domain. The effect of numerical quadrature is also examined.
In the semidiscrete case, optimal error estimates in $L^\infty(L^2)$ and  $L^\infty(H^1)$- 
norms are shown to hold with minimal regularity assumptions on the 
initial data, whereas quasi-optimal estimate in derived in  $L^\infty(L^\infty)$-norm under higher regularity on the data. Based on a second order explicit method in time, a completely discrete scheme is examined and optimal error estimates are established with a mild condition on the space and time discretizing parameters.  Finally, some 
numerical experiments are conducted which confirm the theoretical order of convergence.
}


\maketitle

\section{Introduction}
\noindent In this paper, we discuss and analyze a finite volume element 
method for approximating solutions to  the following class of second order linear hyperbolic integro-differential equations: 
\begin{eqnarray}
\label{a} 
u_{tt}-\nabla\cdot\left(\cA(x)\nabla u+\int_{0}^{t} \cB(x,t,s) \nabla u(s) \,ds \right)&=&f(x,t)
\;\;\;\;\;\mbox{in}~~\Omega\times J, \nonumber \\ 
u(x,t) &=&0\;\;\;\;\;\;\;\;\;\;\;\;\;\mbox{on}~\partial \Omega\times J, \\
u(x,0) &=&u_0(x)\;\;\;\;\;\;\;\mbox{in}~\Omega,  \nonumber\\
u_t(x,0) &=&u_1(x)\;\;\;\;\;\;\;\mbox{in}~\Omega,  \nonumber
\end{eqnarray}
with given functions $u_0$ and $u_1$, where $\Omega\subset \mathbb  {R}^2$ is a bounded
convex polygonal domain, $J=(0,T],~T<\infty$, $u_{tt}=\partial^2 u/
{\partial t^2}$ and $f$ is given function defined on the space-time domain $\Omega\times J.$  Here, $\cA=[{a_{ij}(x)}]$ and 
$\cB=[{b_{ij}(x,t,s)}]$ are $2\times 2$ matrices 
with smooth coefficients. Further, 
assume that $\cA$ is symmetric and uniformly positive definite in $\bar{\Omega}$. 
Problems of this kind arise in linear viscoelastic models, specially in the modelling of viscoelastic materials with memory (cf. Renardy  {\it et al.} \cite{RHN}).

Earlier, the finite volume difference methods which are based on cell centered grids and 
approximating the derivatives by difference quotients have been proposed and analyzed, 
see \cite{EGH}  for a survey. Another approach, which we shall follow in this article 
was formulated in the framework of Petrov-Galerkin finite element method using two different 
grids to define the trial space and test space. This is popularly known finite volume element 
methods (FVEMs). Here and also in literature, 
the trial space consists of $C^0$- piecewise linear polynomials on the finite element partition $\cT_h$ of $\overline{\Omega}$  and the test 
space is  piecewise constants over the control volume $\cT_h^*$ to be defined in
Section 2. Earlier, the FVEM has been examined by Bank and Rose \cite{5}, Cai
\cite{10},
Chatzipantelidis \cite{33}, Li {\it et al.} \cite{19},  Ewing {\it et al.} \cite{15}, etc.  for elliptic problems, 
for parabolic and parabolic type problems  by Chou {\it et al.} \cite{6}, Chatzipantelidis {\it et al.} \cite{25},
Ewing et al. \cite{50}, Sinha {\it et al.} \cite{26} and for second order wave equations by Kumar {\it et al.}
\cite{KNP-2008}. For a recent survey on FVEM, see, a review article by Lin {\it et al.} \cite{LLY}.

For linear elliptic problems,  Li {\it et al.}~\cite{19} have established 
optimal error estimates in $H^1$ and  $L^2$-norms. More precisely, for $L^2$-norm the following estimate
are derived:
$$\|u-u_h\|_{0} \leq Ch^2 \|u\|_{W^{3,p}(\Omega)},~~~~~~p>1,$$
where $u$ is the exact solution and $u_h$ is the finite volume element  approximation of
$u$. Compared to the error analysis of finite element methods, it is observed that  this method 
is optimal in approximation property, but is not optimal with respect to the regularity of the 
exact solution as for $O(h^2)$ order convergence, the exact solution $u\in H^3.$ For convex 
polygonal domain $\Omega$, it may be difficult  
to  prove $H^3$-regularity for the solution  $u.$ Therefore, an attempt has been  made 
in \cite{15} to establish optimal $L^2$ error  estimate  under the assumption  that the exact 
solution $u\in H^2$ and the source term $f\in H^1.$  A counter example has also been provided 
in \cite{15} to show that if  $f\in L^2 $, then FVE solution may not have optimal error estimates 
in $L^2$ norm. The analysis has been extended to parabolic problems in  convex polygonal domain 
in \cite{25} and optimal error estimates have been derived under some compatibility  conditions 
on the initial data. Further, the effect of quadrature, that is, when  the $L^2$ inner product 
is replaced by numerical quadrature  has been analyzed. Subsequently, 
Ewing et al. \cite{50} have  employed FVEM for approximating solutions of parabolic integro-differential equations and derived  optimal  error estimates  under $L^{\infty}(H^3)$ regularity for the exact solution and $L^2(H^3)$ regularity for its time derivative. Then on convex polygonal  domain, Sinha {\it et al.} \cite{26}
have examined  semidiscrete FVEM and  proved optimal error estimates for smooth and non smooth data. 
The analysis is further generalized  to a second order linear wave equation defined on a convex 
polygonal domain and {\it a priori} error estimates have been established only for semidiscrete case, 
see, Kumar {\it et al.} \cite{50}. Further, the effect of quadrature and maximum norm estimates are 
proved under some additional  conditions on the initial data and the forcing function. 
In the present article, an attempt has been made to extend the analysis of FVEM  to  a class of 
second order linear hyperbolic integro-differential equations in convex polygonal domains with 
minimal regularity assumptions on the initial data.  Moreover, a completely discrete scheme based 
on a second order explicit method has been analyzed.

In order to put the present investigation into  a proper perspective visa-vis earlier results, we discuss, below, the literature for the second order hyperbolic
equations. Li et al. \cite{19}  have proved  an optimal order of convergence in $H^1$-norm  without quadrature using elliptic projection, but  the  regularity of the exact
solution assumed to be higher than the regularity assumed in our results even when $B=0$ for the problem (\ref{a}). On a related  finite element analysis for the second order hyperbolic equations without quadrature, we refer to  Baker \cite{37}  and   with quadrature, see, Baker and Dougalis \cite{35} and Dupont \cite{38}.
Baker and Dougalis \cite{35} have proved  optimal  order of convergence 
in $L^\infty( L^2)$ for the semidiscrete finite element scheme, 
provided the initial displacement $u_0\in H^5\cap H_0^1$ and the initial velocity $u_1\in H^4\cap H_0^1.$
Subsequently, Rauch \cite{36} has derived the
convergence analysis for the Galerkin finite element methods when applied to a second  order wave
equation by using piecewise linear polynomials and established optimal
$L^\infty (L^2)$ estimate with $u_0\in H^3\cap H_0^1$ and $u_1=0$ which are  turned out to be  the minimal regularity conditions for the second order wave equation. Subsequently, Pani {\it et al.} \cite{39} 
have examined the effect of numerical quadrature on finite element 
method for  hyperbolic integro-differential equations
with minimal regularity assumptions on the initial data, that is, $u_0\in H^3\cap H_0^1$ and $u_1\in H^2\cap H_0^1$. On a related article on a linear second order wave equation, we refer to Sinha  \cite{34} and on hyperbolic PIDE, see, \cite{CL-89}.
When FVEM is combined with quadrature for approximating solution of (\ref{a}), we have, in this article, 
proved optimal $L^\infty (L^2)$ estimate with minimal regularity assumptions on the initial data.

The organization of the present paper is as follows: 
Section $2$ deals with some notations, weak formulation and the regularity results for the exact solution.
Section $3$ is devoted to the primary and dual meshes for  finite volume element method and semidiscrete 
FVE approximation to the problem (\ref{a}). Section $4$ focuses on {\it a priori} error estimates for 
the semidiscrete FVE approximations and optimal order of convergence in $ L^2$ and $H^1$ norms are 
established under minimal regularity assumptions on the initial data. 
Further,  quasi-optimal order of convergence in maximum norm  has also been derived. 
Section $5$ is on completely discrete scheme which is based on  a second order explicit scheme in time
and {\it a priori} error estimates are established.
Section $6$ deals with the effect of numerical quadrature and the related error estimates are derived again with minimal regularity assumption on the initial data. Finally in Section $7$, some numerical experiments are conducted which confirm our theoretical order of convergence.

Through out this paper, $C$ is  a
generic positive constant independent of discretising  parameters $h$ and $k.$

\section{Notation and Preliminaries.}
\se
This section is devoted to some notations and  preliminary results related to the weak solution of (\ref{a}).

Let $W^{m,p}(\Omega)$ denote the standard 
Sobolev space with the norm
$$\|u\|_{m,p,\Omega}=\left(\sum_{|\alpha|\leq m
}\|D^{\alpha}u\|_{L^p(\Omega)}^p\right)^{1/p}~~~~\mbox{for} ~1\leq p<\infty,$$
and for $p=\infty$,
$$\|u\|_{m,\infty,\Omega}=\sup_{|\alpha|\leq m}\|D^\alpha u\|_{L^\infty(\Omega)}.~~~~~~~~~~~~~~~~~~$$
When there is no confusion, we denote $\|u\|_{m,p,\Omega}$ by $\|u\|_{m,p}.$
For $p=2$, we simply write $W^{m,2}(\Omega)$ as $H^m(\Omega)$ and denote its norm by
$\|\cdot \|_{m}.$
For a Banach space $X$ with norm $\|\cdot\|_X$ and $ 1\leq p \leq \infty,$ let $W^{m,p}(0,T; X)$  be  defined by 
$$ W^{m,p}(0,T; X):=\{v:(0,T) \longrightarrow X | \|D_t^j v\|_X \in L^p(0,T),\;\; 0\leq j \leq m\}.$$
with its norm
$$
\|v\|_{W{m,p} (0,T; X)}=\|u\|_{W^{m,p}(X)}:=\sum_{j=0}^{m}
\left(\int_0^T\|D^j_t v\|_{X}^p\;dt\right)^{1/p},
$$ 
with the standard modification for $p=\infty$, see \cite{42}. For $m=0$, $W^{m,p}(0,T; X)$ is simply 
the space $L^p(X)$.
Finally, let $(\cdot ,\cdot)$ and $\|\cdot\|_0$ denote, respectively, the $L^2$ inner product and its
induced norm on $L^2(\Omega).$ 

With
$H_0^1(\Omega)=\{v\in H^1(\Omega): v=0 \;\mbox{on}\;  \partial \Omega \},$
define  the bilinear forms $A(\cdot,\cdot)$ and $B(\cdot,\cdot)=B(t,s;\cdot,\cdot)$ 
on $H_0^1(\Omega)\times H_0^1(\Omega)$  by
$$ A(u,v)= \int_{\Omega} \cA(x)\nabla u\cdot \nabla v\,dx,$$
and $$ B(t,s;u(s),v)= \int_{\Omega} \cB(x,t,s)\nabla u(s)\cdot \nabla v\,dx.$$
Then, the weak formulation for (\ref{a}) is to seek  
$u:(0,T]\longrightarrow H^1_0(\Omega)$ such that
\begin{equation} \label{b}
(u_{tt},v)+ A(u,v)+\int_{0}^{t} B(t,s;u(s),v) \,ds =  (f,v)\;\;\;
\forall v\in H_0^1(\Omega)
\end{equation}
with  $u(0)=u_0$ and $u_{t}(0)=u_1.$

Since $\cA$ is symmetric and uniformly positive definite in $\Omega$, 
the bilinear form  $ A(\cdot, \cdot)$ satisfies the following  condition: 
there exist positive  constants $\alpha $ and $\Lambda$ with $\Lambda \geq \alpha$  such that 
\begin{equation}
\Lambda \|v\|_1^2 \geq A(v,v) \geq \alpha \|v\|_1^2\;\;\;\forall v\in H_0^1(\Omega).
\label{eqn2.3}
\end{equation}

For our subsequent use, we state without proof {\it a priori} estimates of 
the solution $u$ of the problem (\ref{a}) under appropriate regularity conditions and compatibility conditions
on $u_0$, $u_1$ and $f$. Its proof can be easily obtained by appropriately  modified arguments in the proof of Theorem 3.1 of \cite{39}.  For similar estimates for second order  linear hyperbolic equations, see  Lemma 2.1 of \cite{KNP-2008}.
\begin{lemma} \label{lem1}
Let $u$ be a weak solution  of $(\ref{a})$. Then, there is a positive constant $C=C(T)$ such that  the following estimates 
\begin{eqnarray*}
\|D^{j+2}_tu(t)\|_0+\|D^{j+1}_t u(t)\|_1 +\|D^j_tu(t)\|_2
& \leq & C \Big( \|u_0\|_{j+2}+\|u_1\|_{j+1}+\\
&& \sum_{k=0}^{j}\displaystyle \|D_t^k f\|_{L^1(H^{j-k})}+ \|D_t^{j+1} f\|_{L^1(L^2)}\Big),
\end{eqnarray*}
hold for $j=0,1,2,$  where $D^j_t=(\partial^j/\partial t^j).$
\end{lemma}
We shall have occasion to use the following identity for $\phi \in C^1 ([0,T];X),$ where 
$X$ is a Banach space
\begin{equation}\label{phi-t}
\phi(t) = \phi(0) + \int_{0}^t \phi_t(s)\; ds.
\end{equation}

\section { Finite Volume Element Method} 
\se
This section deals with primary and dual meshes on the domain $\Omega$, construction of finite dimensional spaces, finite volume element formulation and some preliminary results.  

\begin{figure}
    \begin{center}
    \includegraphics*[width=11.0cm,height=7.0cm]{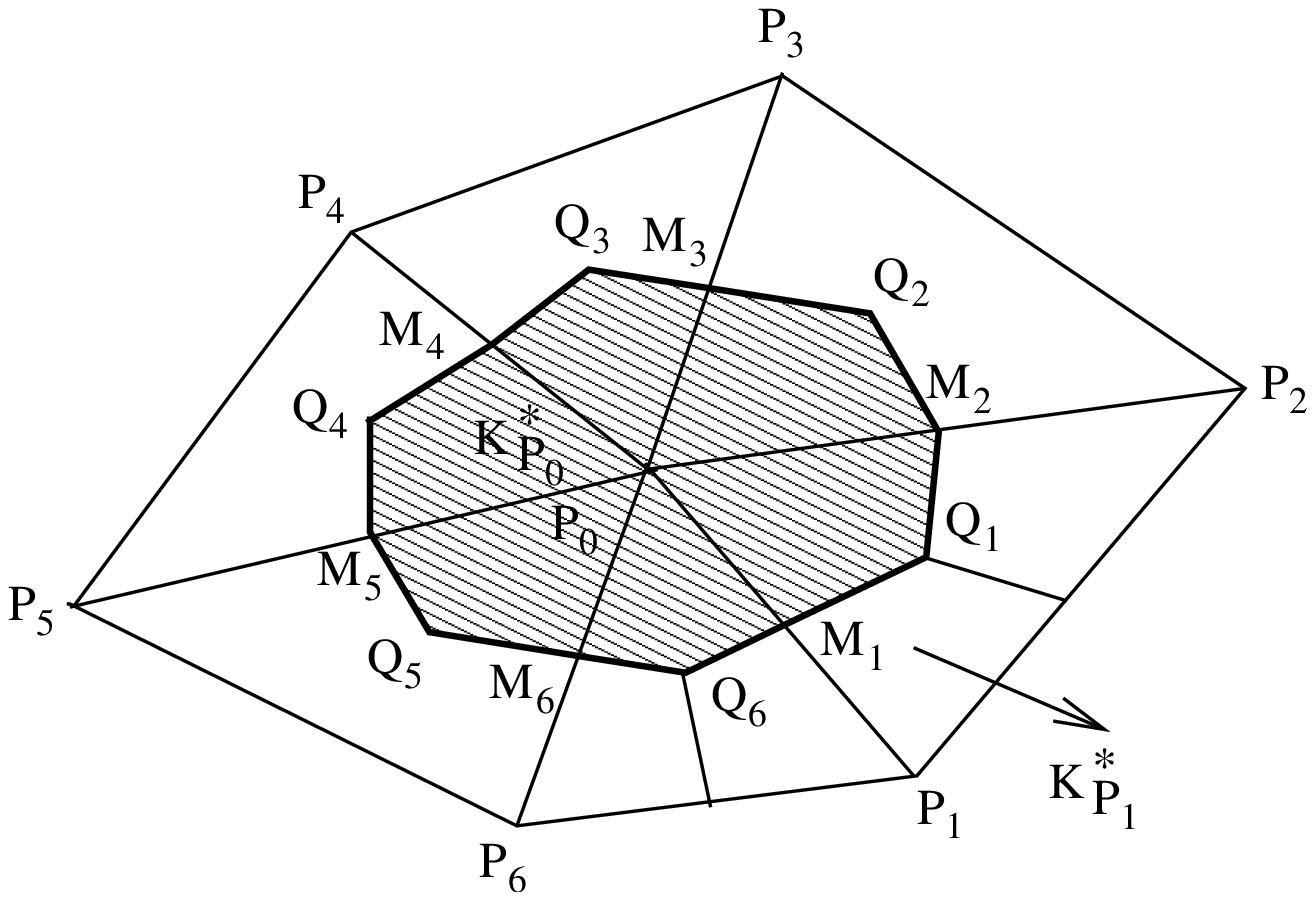}
       \caption{}
       \label{fig:mesh}
    \end{center}
\end{figure}

Let $\cT_h$ be a family of regular triangulations of 
the closed, convex  polygonal domain $\overline{\Omega}$ into closed triangles  $K,$ 
and let $h=\max_{K\in \cT_h}(\mbox{diam}K),$ where $h_{K}$ denotes the diameter  of $K.$  Let $N_h$ be set of nodes or vertices, that is, $N_h :=\left\{P_i:P_i~~\mbox{ is a vertex of the element }~K \in 
\cT_h~\mbox{and}~P_i\in \overline{\Omega}\right\}$ and let
$N_h^0$ be the set of interior nodes in $\cT_h$ with cardinality $N$.
Further, let $\cT_h^*$ be the dual mesh associated with the primary mesh $\cT_h,$ which is defined as follows. With $P_0$ as an interior node of the triangulation $\cT_h,$ let  $P_i\;(i=1,2\cdots
m)$ be its adjacent nodes (see, FIGURE~\ref{fig:mesh} with $m=6$ ). Let $M_i,~i=1,2\cdots
m$ denote the midpoints of $\overline{P_0P_i}$ and let $Q_i,~i=1,2\cdots
m,$  be  the barycenters of the triangle $\triangle P_0P_iP_{i+1}$ with
$P_{m+1}=P_1$. The {\it control volume}
  $K_{P_0}^*$ is constructed  by joining successively $ M_1,~ Q_1,\cdots
  ,~ M_m,~ Q_m,~ M_1$. With $Q_i ~(i=1,2\cdots
m)$ as the nodes of $control~volume~$ $K^*_{p_i},$ let $N_h^*$ be the set of all dual nodes
$Q_i$. For a  boundary
node $P_1$,  the control volume $K_{P_1}^*$ is shown in the FIGURE~\ref{fig:mesh}. Note that the union
of the control volumes forms a partition $\cT_h^*$ of $\overline{\Omega}$.

Assume that the partitions $\cT_h$ and $\cT_h^*$ are
quasi-uniform in the sense that  there exist positive constants $C_1$ and $C_2$
independent of $h$ such that 
\begin{eqnarray}
C_1\;h^2\leq |K_{Q_i}| \leq C_2\; h^2~~~~~\forall Q_i \in N_h^*,
 \label{p1}
 \end{eqnarray}
 \begin{eqnarray}
 C_1\;h^2\leq |K^*_{P_i}|\leq C_2 \;h^2~~~~~\forall P_i \in N_h,
\end{eqnarray}
where $ |K|= \mbox {\; meas\;} (K).$

We consider a finite volume element discretization of (\ref{a}) in the standard $C^0$-conforming piecewise linear  finite element space $U_h$ on the primary mesh $\cT_h$, which is defined by
$$U_h=\{v_h\in C^0(\overline {\Omega})\;:\;v_h|_{K}\;\mbox{is linear for all}~ K\in \cT_h\; \mbox{and} \; v_h|_{\partial \Omega}=0\},$$
and  the dual volume element space $\vv$ on the dual mesh $\cT^*_h$ given by 
$$\vv=\{v_h\in L^2(\Omega)\;:\;v_h|_{K_{P_0}^*}\;\mbox{is constant for all}\; K_{P_0}^*\in \cT_h^*\; \mbox{and}\; v_h|_{\partial \Omega}=0\}.$$
Now, $U_h=\mbox{span} \{\phi_i\;:\;P_i\in N_h^0\}$ and 
$\vv=\mbox{span}\{\chi_i\;:\; P_i\in N_h^0\}$, where $\phi_i$'s are the standard nodal 
basis  functions associated with nodes $P_i$ and 
$\chi_i$'s are the characteristic basis functions corresponding  to the control volume 
$K_{P_i}^*$ given by 
$$\chi_i(x)=\left\{\begin{array}{cc}
1,&\mbox{if}~ x\in K_{P_i}^*\\
0,& \mbox{elsewhere}.\\
\end{array}\right.$$

The semidiscrete finite volume element formulation for (\ref{a}) is  
to seek  $u_h:(0,T]\longrightarrow U_h$ such that 
\begin{equation}
\label{c} 
(u_{h,tt},v_h)+ A_h(u_h,v_h)+\int_{0}^{t} B_h(t,s;u_h(s),v_h) \,ds=  (f,v_h)\;\;
\forall v_h\in \vv,~~~~~~~~~~~~~~~~~~~~~ 
\end{equation}
with given initial data $u_h(0)$ and $u_{h,t}(0)$ in $U_h$ to be defined later.
Here, the bilinear forms $A_h(\cdot, \cdot)$ and $B_h(t,s;\cdot, \cdot)$ 
are  defined, respectively, by
$$ A_h(u_h,v_h)= -\sum_{P_i\in N_h^0} v_h(P_i) \int_{\partial K_{P_i}^*}
\cA(x) \nabla u_h\cdot{\bf n}\,ds,~~~~$$
and 
$$B_h(t,s;u_h,v_h)= -\sum_{P_i\in N_h^0} v_h(P_i) \int_{\partial K_{P_i}^*}
\cB(x,t,s) \nabla u_h\cdot{\bf n}\,ds$$
for all $(u_h,v_h)\in U_h \times \vv$, 
 with ${\bf n}$ denoting the outward unit normal to the boundary of the control volume 
$K_{P_i}^*$. Notice that by taking the $L^2$ inner product of (\ref{a}) with 
$v_h\in \vv$ and then integrating, we obtain a similar equation for $u$ as 
\begin{equation}\label{newaa1}
(u_{tt},v_h)+ A_h(u,v_h)+\int_{0}^{t} B_h(t,s;u(s),v_h) \,ds =  
(f,v_h)\;\;\;\forall v_h \in \vv.
\end{equation}

For the error analysis, we first introduce two interpolation operators.
Let $\Pi_h:C(\Omega)\longrightarrow U_h$ be the piecewise linear interpolation operator and $\Pi_h^*:C(\Omega)\longrightarrow \vv$ be the piecewise constant interpolation operator. These interpolation operators are defined,  respectively,  by
\begin{equation}
\Pi_h u=\sum_{P_i\in N_h^0}u(P_i)\phi_i(x)\;~\mbox{and}~\; \Pi_h^*u=\sum_{P_i\in N_h^0}u(P_i)\chi_i(x).
\label{naa}
\end{equation}
Now for $\psi \in H^2$,  $\Pi_h$ has the following approximation property, (see, 
Ciarlet \cite{46}):
\begin{equation}
\|\psi-\Pi_h\psi\|_0\leq Ch^2\|\psi\|_2.
\label{1.6}
\end{equation}
Further, we introduce the following discrete norms
$$\|v_h\|_{0,h}=\left(\sum_{K\in T_h}|v_h|_{0,h,K}^2\right)^{1/2} \;\; \mbox {and}\;\;	 \|v_h\|_{1,h}=\left(\|v_h\|_{0,h}^2+|v_h|_{1,h}^2\right)^{1/2},
$$
where the seminorm 
$|v_h|_{1,h}=\left(\sum_{K\in T_h}|v_h|_{1,h,K}^2\right)^{1/2}$, and for  $K=K_Q=\triangle P_1P_2P_3$,
$$|v_h|_{0,h,K}=\left\{\frac {1}{3}\left(v_h(P_1)^2+v_h(P_2)^2+v_h(P_3)^2\right)\;|K|\;\right\}^{1/2}
$$
$$
|v_h|_{1,h,K}=\left\{(|\frac{\partial v_h}{\partial x}|^2+|\frac{\partial v_h}{\partial y}|^2
)\; |K|\;\right\}^{1/2}.
$$
In the following Lemma, a relation between discrete norms and  standard Sobolev  norms
is stated without  proof. For a proof, see,  \cite[pp. 124]{19} and \cite{10}.
\begin{lemma}
\label{lemma1}
For $v_h \in U_h,~ |\cdot|_{1,h}$ and $|\cdot|_1$ are identical; $\|\cdot\|_{0,h}$ and $\|\cdot\|_{1,h}$ are equivalent to  $\|\cdot\|_{0}$ and $\|\cdot\|_{1}$,  respectively, that is, there exist  positive constants $C_3$ and $C_4>0$, independent of h, such that
\begin{eqnarray}
C_3\|v_h\|_{0,h}\leq \|v_h\|_0\leq C_4\|v_h\|_{0,h}\;\;\; \forall v_h\in U_h
\end{eqnarray}
and
\begin{equation}
C_3||v_h||_{1,h}\leq ||v_h||_1\leq C_4||v_h||_{1,h}\;\;\;\;\forall v_h\in U_h.
\label{nb}
\end{equation}
\end{lemma}
Note that $ \|v_h\|_{0,h} = \|\Pi_h^* v_h\|_0.$  Below, we state without proof the properties of  
the interpolation operator $\Pi_h^*.$ For a proof, we refer to  \cite[pp. 192]{19}.
\begin{lemma}
\label{E}
 The following statements  hold true.\\
 $ (i) \; \mbox{For} \;\;\Pi_h^*:U_h\longrightarrow \vv  
 \mbox{ defined in (\ref{naa}),} \qquad \qquad \qquad \qquad \qquad$
\begin{equation}
  (\phi_h,\Pi_h^*v_h)=(v_h,\Pi_h^*\phi_h)\;\;\;\;\;\forall \phi_h,~v_h\in
 U_h.
\label{nc}
\end{equation}
(ii) \;With $\||\phi_h\||:= (\phi_h,\Pi_h^*\phi_h)^{1/2}$, the norms  $\||\cdot\||$  and  
$\|\cdot\|_0 $  are equivalent on $ U_h$,  that is, there exist positive constants $c_{eq}$ and  $C_{eq}$, independent of $ h$,  such that
\begin{equation}
c_{eq}||\phi_h||_{0}\leq \||\phi_h\||\leq C_{eq}||\phi_h||_{0}\;\;\;\;\;\;\;\forall \phi_h\in U_h.
\label{nd}
\end{equation}
 \end{lemma} 
 
\section{A Priori Error Estimates}
\se
This section is devoted to {\it a priori} error estimates of the approximation $u_h$ to the spatial semidiscrete scheme (\ref{c}). 

For the derivation of optimal error estimates, we split  $e=u-u_h$ as 
$$e:=(u-\rv u)+(\rv u-u_h)=:\rho+\theta,$$  
where  $\rv:L^{\infty}(H_0^1\cap H^2)\rightarrow L^{\infty}(U_h)$ is the Ritz-Volterra projection 
defined by
\begin{equation}\label{d}
A(u-\rv u,\chi_h)+\int_{0}^{t} B(t,s;u-\rv u,\chi_h)\,ds=0\;\;\;\;\forall \chi_h\in U_h.
\end{equation}
With some abuse of notations, we will denote by $\rv u_0$ the Ritz  projection of $u_0$ onto $U_h$ defined by
$$
A(u_0-\rv u_0,\chi_h)=0\;\;\;\;\forall \chi_h\in U_h.
$$
For our subsequent analysis, we state without proof following error estimates for the Ritz-Volterra projection.
For a proof, see, \cite{39}, \cite{CL-88}, \cite{YL}, \cite{LTW} and \cite{PTW}.
\begin{lemma}  
\label{lem2}
There exist  positive constants $C$, independent of $h$,  such that
for $j=0,1,2$, and $r=1,2$ the following estimates hold:
\begin{equation}\label{1Rh}
\| D_t^j\rho(t)\|_0+h\| D_t^j\rho(t)\|_1\leq C h^r
\left[\sum_{l=0}^j\| D_t^lu(t)\|_r+ \int_0^t\|u(s)\|_r\;ds\right],
\end{equation}
and 
\begin{equation}\label{R1h}
\|\rho(t)\|_{0,\infty}\leq Ch^2 \left(\log\frac{1}{h}\right)\left(\|u(t)\|_{2,\infty}
+\int_0^t||u(s)||_{2,\infty}ds\right).
\end{equation}
\end{lemma}

Now, define 
$$\epsilon_h(f,\chi)=(f,\chi)-(f,\Pi_h^*\chi)\;\;\;\;\forall \chi
\in U_h,~~~~~~~~~~~~~~~$$
$$\;\;\;\;\;\;\;\;\;\;\;\;\;\epsilon_A(\psi,\chi)=A(\psi,\chi)-A_h(\psi,\Pi_h^*\chi)\;\;\;\;\forall \psi,~ \chi
\in U_h,$$
and 
$$\;\;\;\;\;\;\;\;\;\;\;\;\;\epsilon_B(t,s;\psi,\chi)=B(t,s;\psi,\chi)-
B_h(t,s;\psi,\Pi_h^*\chi)\;\;\;\;\forall \psi,~ \chi
\in U_h.$$ 
Then, the following lemma will be of frequent use in our analysis and the proof of which 
can be found in \cite{33}.
\begin{lemma}\label{A}
Assume that the coefficient matrices $\cA, \cB(t,s) \in W^{1+i,\infty} (\Omega; \mathbb{R}^{2\times 2})$  for $i=0,1.$ Then, there exist positive constant $C,$ independent of $h$, such that the following estimates hold for
$\chi \in U_h$  and for $i,~j=0,~1$
\begin{eqnarray} \label{f}
 |\epsilon_h(f,\chi)|\leq Ch^{i+j}\|f\|_{H^i}\;\|\chi \|_{H^j}\;\;\;
 \forall f\in H^i,
\label{ne}
\end{eqnarray}
and for $u\in H^{1+i}\cap H_0^1$
\begin{eqnarray}\label{ea-1}
 |\epsilon_A(\rv u,\chi)|  \leq Ch^{i+j}\Big(\|u\|_{H^{1+i}} + \int_{0}^{t} \|u(s)\|_{H^{1+i}}\;ds\Big)\;\|\chi \|_{H^j}.
\label{nf}
\end{eqnarray}
Moreover,  
\begin{eqnarray}\label{ea-2}
|\epsilon_A(w_h,\chi)|  \leq Ch\|w_h\|_{H^{1}}\;\|\chi
\|_{H^1}\;\forall w_h \in U_h.
\label{ng}
\end{eqnarray}
\end{lemma}
\noindent
The estimates (\ref{ea-1}) and (\ref{ea-2}) are also valid if $\epsilon_A$ is replaced by $\epsilon_B.$

Now, for $\psi\in H^1_0$ and for each $t\in (0,T],$ introduce a linear functional $\G(\psi)=\G(t,\psi)$ defined on $U_h$ by
$$ \G(\psi)(\chi)=\epsilon_A(\psi,\chi)+\int_0^t\epsilon_B(t,s;\psi(s),\chi)\,ds,\,\,\,\chi \in U_h.$$
Notice that, by using the definition of $\G$, (\ref{b}) and  (\ref{newaa1}), 
there follows that
\begin{eqnarray}
\G(\rho)(\chi)
&=&A(u,\chi)+\int_{0}^{t} B(t,s;u(s),\chi)\,ds\nonumber \\
&&-A_h(u,\Pi_h^*\chi)-\int_{0}^{t} B_h(t,s;u(s),\Pi_h^*\chi)\,ds
-\G(\rv u)(\chi)\nonumber \\
&=&(f-u_{tt},\chi)-(f-u_{tt},\Pi_h^*\chi)-\G(\rv u)(\chi)\nonumber \\
&=&\epsilon_h(f-u_{tt},\chi)-\G(\rv u)(\chi).
\label{h}
\end{eqnarray}
From (\ref{c}) and (\ref{newaa1}),  we obtain the equation in $\theta$ for $v_h\in \vv$ as
$$
(\theta_{tt},v_h)+A_h(\theta,v_h)+\int_{0}^{t} B_h(t,s;\theta(s),v_h)\,ds=
-A_h(\rho,v_h)-\int_{0}^{t} B_h(t,s;\rho,\chi)\,ds-(\rho_{tt},v_h).
$$
Choosing $v_h=\Pi_h^*\chi $ and using the definition of $\G$ 
and  (\ref{d}), we find  that 
\begin{eqnarray} \label{i}
(\theta_{tt},\Pi_h^*\chi)+A(\theta,\chi)\;ds&+&\int_{0}^{t} B(t,s;\theta(s),\chi) \,ds = G(\rho)(\chi)\nonumber\\
&+& G(\theta)(\chi)-(\rho_{tt},\Pi_h^*\chi) \;\;\forall \chi \in U_h.
\end{eqnarray}

For any continuous function $\phi$ in $[0,t]$,  define 
$\hat{\phi}$ by 
$$ \hat{\phi}(t)=\int_0^t\phi(s)\,ds.$$
Notice that $\hat{\phi}(0)=0$ and $(d\hat{\phi}/dt)(t)=\phi(t)$.
Then, integrate  (\ref{i}) from $0$ to $t$ to obtain 
\begin{eqnarray}\label{ii}
(\theta_{t},\Pi_h^*\chi)+A(\hh,\chi)&=&
 \bG(\rho)(\chi)+\bG(\theta)(\chi)
+(-\rho_{t},\Pi_h^*\chi)+(e_t(0),\Pi_h^*\chi)\nonumber\\
&&-\int_{0}^{t} B(s,s;\hh(s),\chi)ds+
\int_{0}^{t}\int_{0}^{s}B_\tau(s,\tau;\hh(\tau),\chi)d\tau ds,
\end{eqnarray}
where 
$$
\bG(\phi)(\chi)=\epsilon_A(\hat{\phi},\chi)+
\int_{0}^{t} \epsilon_B(s,s;\hat{\phi}(s),\chi)ds-\int_{0}^{t}
\int_{0}^{s}\epsilon_{B_\tau}(s,\tau;\hat{\phi}(\tau),\chi)d\tau ds.
$$
For a linear functional $F$ defined on $U_h$, set
$$\|F\|_{-1,h}=\sup_{0\neq \chi\in U_h}\frac{|F(\chi)|}{\|\chi\|_1}.$$
We shall need the following lemmas in our subsequent analysis. 
\begin{lemma}\label{lm-n}
With $\G$ and $\bG$ as above, there exists a positive constant $C=C(T)$ such that the following estimates
\begin{equation} \label{G-1}
\|D_t^j\G(\rv u) \|_{-1,h}\leq Ch^2  \left( \sum_{\ell=0}^{j} \|D^{\ell}_t u(t) \|_{2}+ \int_{0}^t \| u (s)\|_{2}\,ds\right),
\end{equation}
and
\begin{equation} \label{G-2}
\|D_t^j\bG(\rv u) \|_{-1,h}\leq C h^2
\left(\sum_{\ell=0}^{j} \|D^j_t \hat {u} (t)\|_{2}+ \int_{0}^{t} \|\hat{u}(s)\|_{2}\; ds\right),
\end{equation}
hold for $j=0,1$.
\end{lemma}
\noindent Proof. Using (\ref{nf}) and the estimates in Lemma~\ref{lem1}, we obtain
\begin{eqnarray*}
|\G(\rv u)(\chi)|&\leq&|\epsilon_A(\rv u,\chi)|+\int_0^t|\epsilon_B(t,s;\rv u(s),\chi)|\,ds\\
&\leq&Ch^2\left(\|u\|_2+\int_0^t\|u(s)\|_2\,ds\right)\|\chi\|_1,
\end{eqnarray*}
and
\begin{eqnarray*}
|\G_t(\rv u)(\chi)|
&\leq&Ch^2\left(\|u_t\|_2+\|u\|_2 + \int_0^t\|u(s)\|_2\,ds\right)\|\chi\|_1\\
&\leq&Ch^2 \left(\|u_t\|_2+ \|u\|_2 + \int_0^t\|u(s)\|_2\,ds\right)\|\chi\|_1
\end{eqnarray*}
In a similar manner,  we derive the second estimate (\ref{G-2}) and this completes the rest of the proof.
\hfill{{\rule{2.5mm}{2.5mm}}}

In the error analysis, we shall frequently use the following inverse assumption:
\begin{equation}\label{inv}
\|\chi\|_1\leq C_{inv} h^{-1}\|\chi\|_0,\quad \chi \in U_h.
\end{equation}

\subsection{ $H^1$-  error estimate}
\begin{theorem}
\label{H2}
Let $u$ and $u_h$  be the solutions of $(\ref{a})$ and $(\ref{c}),$ respectively, and assume 
that $f\in L^1(H^1),\; f_t,\;f_{tt}\in L^1(L^2)$, $u_0\in H^3\cap H_0^1$ and $u_1\in H^2\cap H_0^1$. Further,
assume that  $u_h(0)=\Pi_hu_0$ and $u_{h,t}(0)=\Pi_hu_1$, where $\Pi_h$ is the
interpolation operator defined in $(\ref{naa})$. Then, there exists a positive constant $C=C(T)$, independent of $h$, such that for $t\in (0,T]$ the following estimate 
\begin{eqnarray*}
 \|u(t)-u_h(t)\|_1 \leq C\;h \left(\|u_0\|_3+\|u_1\|_2+ \int_{0}^{t} \Big( \|f\|_{1} +\|f_t\|_0 + \|f_{tt}\|_0\Big)\,ds \right)  
 \end{eqnarray*}
holds.
\end{theorem}
\noindent Proof. Since $u-u_h=\rho+\theta$ and estimates of $\rho$ are known from the Lemma~\ref{lem2}, 
it is sufficient to estimate $\theta.$  Choose  $\chi=\theta_t$ in (\ref{i}) and use (\ref{h}) to obtain   
\begin{eqnarray}
(\theta_{tt},\Pi_h^*\theta_t)+A(\theta,\theta_t)+\int_{0}^{t} B(t,s;\theta(s),\theta_t) \,ds
&=& \epsilon_h(f-u_{tt},\theta_t)
-\G(\rv u)(\theta_t) \nonumber\\
&&+\G(\theta)(\theta_t)-(\rho_{tt},\Pi_h^*\theta_t).\nonumber
\label{k}
\end{eqnarray}
Now use (\ref{nc}) and  symmetry of the  bilinear form $A(\cdot,\cdot)$ to arrive at    
\begin{eqnarray*}
\frac{1}{2}\frac{d}{dt}\Big[(\theta_t,\Pi_h^*\theta_t)+A(\theta,\theta)\Big] &=& \epsilon_h(f-u_{tt},\theta_t)
-\G(\rv u)(\theta_t)+\G(\theta)(\theta_t)-(\rho_{tt},\Pi_h^*\theta_t)\\
&&-\int_{0}^{t} B(t,s;\theta(s),\theta_t(t)) \,ds.
\end{eqnarray*}
Integration from $0$ to $t$ yields
\begin{eqnarray}
\frac{1}{2}\Big(\||\theta_t\||^2+A(\theta,\theta)\Big)&=& \frac{1}{2}\left(\||\theta_t(0)\||^2+A(\theta(0),\theta(0))\right)+
\int_0^t \epsilon_h(f-u_{tt},\theta_t)ds \nonumber \\
&&- \int_0^t \G(\rv u)(\theta_t) ds+ 
 \int_0^t \G(\theta)(\theta_t)ds + \int_0^t(-\rho_{tt},\Pi_h^*\theta_t)ds\nonumber \\
&&-\int_{0}^{t} \int_0^s B(s,\tau;\theta(\tau),\theta_t(s)) \,d\tau ds\nonumber\\
&=& J_1+J_2+J_3+J_4+J_5+J_6.
\label{1.27}
\end{eqnarray}
For the first term on the right hand side of (\ref{1.27}), a use of the boundedness 
of $A(\cdot,\cdot)$ with (\ref{1.6}) and  (\ref{nd}) shows 
\begin{equation}
|J_1|\leq C \left(\|\theta_t(0)\|_0^2+\|\theta(0)\|_1^2\right)\leq Ch^2(\|u_1\|_1^2+\|u_0\|_2^2).
\label{1.28}
\end{equation}
For estimating $J_2,$  an application of  (\ref{ne}) with $j=0$  
implies 
\begin{eqnarray}
|J_2| \leq Ch 
\int_0^t\left(\|f\|_1+\|u_{tt}\|_1\right)\|\theta_t\|_0\;ds.
\end{eqnarray}
To estimate  $J_3$, a use of the inverse inequality (\ref{inv}) shows that
\begin{equation}
|J_3|\leq C \int_{0}^{t}\|\G(\rv u)\|_{-1,h}\|\theta_t\|_1 \;ds \leq 
Ch^{-1}\int_{0}^{t}\|\G(\rv u)\|_{-1,h}\|\theta_t\|_0\;ds.
\end{equation}
Using the definition of $G$, (\ref{ng}) and  the inverse inequality, it follows that 
\begin{eqnarray}
\;\;\;\;\;\;\;\;\;|J_4|&\leq &
\int_0^t|\G(\theta)(\theta_t)| \;ds \nonumber\\
&\leq& Ch\left[\int_0^t\|\theta\|_1\|\theta_t\|_1 ds +\int_0^t\int_0^s\|\theta(\tau)\|_1
\|\theta_t(s)\|_1d\tau ds\right] \nonumber\\
&\leq& C\left[\int_0^t\|\theta\|_1\|\theta_t\|_0 ds+\left(\int_0^t\|\theta\|_1 ds\right)
\left(\int_0^t\|\theta_t\|_0 ds\right)\right]. 
\end{eqnarray}
For $J_5$, apply the Cauchy-Schwarz inequality, $L^2$ stability of
$\Pi_h^*$ and (\ref{1Rh}) with $r=1$ to obtain 
\begin{equation}
|J_5|\leq  \int_{0}^{t}\|\rho_{tt}\|_0\;\|\theta_t\|_0 \;ds
\leq C(T)\;h\int_{0}^{t}\Big(\|u_{tt}\|_1+ \|u_{t}\|_1+\|u\|_1\Big)\|\theta_t\|_0\; ds.
\label{1.29}
\end{equation}
For the term $J_6,$ we note that an integration by parts  yields 
\begin{eqnarray*} \label{integral-term}
\int_0^t \int_0^s B(s,\tau;\theta(\tau),\theta_t(s))\,d\tau \;ds &=&
\int_0^tB(t,s;\theta(s),\theta(t))\,ds - \int_0^tB(s,s;\theta(s),\theta(s)) \\
&&-\int_0^t\int_0^sB_s(s,\tau;\theta(\tau),\theta(s))\,d\tau \;ds, \nonumber
\end{eqnarray*}
and hence, deduce that
\begin{equation}
|J_6|\leq C\left(\|\theta(t))\|_1\;\int_0^t\|\theta(s)\|_1\,\;ds
+ \int_0^t\|\theta(s)\|^2_1\,ds\right).
\label{eer}
\end{equation}
Now, set $\cE^2_1(t)=\|\theta_t\|^2_0+\|\theta\|^2_1$ and
$$\cE_1(t^\ast)=\max_{0\leq \tau\leq t}\cE_1(t),$$ for some $t^\ast\in [0,t]$. Then, 
substituting the estimates  (\ref{1.28})-(\ref{eer}) in (\ref{1.27}),   
using coercivity of $A(\cdot,\cdot)$,  equivalence of norms $\||\cdot\||$ and $\|\cdot \|_0,$ 
apply  standard kick back arguments to find that 
\begin{eqnarray*}
\cE_1(t^\ast) &\leq & Ch\;\Big( \|u_0\|_2+\|u_1\|_1+\int_{0}^{T} \left(\|u_{tt}(s)\|_1 + \|u_{t}(s)\|_1 
+\|u(s)\|_1\right)\;ds\Big)\\
 &+& C h \;\int_{0}^{T} \left( \|f(s)\|_1 + h^{-2}\|\G(\rv u)(s)\|_{-1,h}\right)\;ds
 + \int_{0}^{t^\ast} \cE_1(s)\; ds.
 \end{eqnarray*}
Now replace $t^\ast$ by $t$ and apply Gronwall's lemma with the estimate (\ref{G-1}) to conclude
that
$$
\cE_1(t) \leq~Ch\Big(\|u_0\|_2+\|u_1\|_1+ \int_{0}^{T}
 \left(\|u\|_2 +\|u_t\|_1 +\|u_{tt}\|_1+\|f\|_1\right)\;ds\Big).
$$
 A use of  triangle inequality with (\ref{1Rh}) and the estimates from Lemma~\ref{lem1} completes 
 the rest of the proof. \hfill{{\rule{2.5mm}{2.5mm}}}


\subsection{ Optimal $L^2$- error estimates }
In this subsection, we shall discuss optimal $L^{\infty}(L^2)$ estimates
\begin{theorem}
\label{TH1}
Under the assumptions of Theorem~$\ref{H2},$ 
 there exists a positive constant $C=C(T)$, independent of $h$, such that 
$$
\|u(t)-u_h(t)\|_0\leq C h^2 \left(\|u_0\|_3+\|u_1\|_2 + \int_{0}^{t} \Big( \|f\|_{1} +\|f_t\|_0 + \|f_{tt}\|_0\Big)\,ds \right).
$$
\end{theorem}
\noindent Proof. 
By setting $\chi=\theta$ in (\ref{ii}), and
using (\ref{nc}) with symmetry of the bilinear form $A(\cdot,\cdot)$, we find that
\begin{eqnarray*}
\frac{1}{2}\frac{d}{dt}\left[(\theta,\Pi_h^*\theta)+A(\hat{\theta},\hat{\theta})\right]&=& 
 \bG(\rho)(\theta)+\bG(\theta)(\theta)
+(-\rho_{t},\Pi_h^*\theta)+(u_1-\Pi_hu_1,\Pi_h^*\theta)\nonumber\\
&&-\int_{0}^{t} B(s,s;\hh(s),\theta(t))\;ds+
\int_{0}^{t}\int_{0}^{s}B_\tau(s,\tau;\hh(\tau),\theta(t))\;d\tau \;ds.
\end{eqnarray*}
Integrate from $0$ to $t$ to obtain 
\begin{eqnarray}\label{nn1}
\frac{1}{2}\left[ \||\theta\||^2+A(\hat{\theta},\hat{\theta})\right] &=&  
\frac{1}{2}\||\theta(0)\||^2 + 
 \int_0^t\bG(\rho)(\theta)\;ds+\int_0^t\bG(\theta)(\theta)\;ds
+\int_0^t(-\rho_{t},\Pi_h^*\theta)\;ds\nonumber\\
&&+(u_1-\Pi_hu_1,\Pi_h^*\hat{\theta})
-\int_{0}^{t}\int_{0}^{s} B(\tau,\tau;\hh(\tau),\theta(s))\; d\tau \;ds\nonumber\\
&&+\int_{0}^{t}\int_{0}^{s}\int_{0}^{\tau'}B_{\tau'}(\tau,\tau';\hh(\tau'),\theta(s))\;d\tau'
\;d\tau \;ds\nonumber \\
&=&\frac{1}{2}\||\theta(0)\||_0^2+I_1+I_2+I_3+ I_4+I_5+ I_6.
\end{eqnarray}
To estimate $I_1$, we note  from (\ref{h}) that
\begin{eqnarray}\label{keyy1}
\bG(\rho)(\theta)
&=&\epsilon_h(\hat{f}-\hat{u}_{tt},\theta)-\bG(\rv u)(\theta)\nonumber\\
&=&\frac{d}{dt}\left (\epsilon_h(\hat{f}-\hat{u}_{tt},\hat{\theta})-\bG(\rv u)(\hat{\theta})\right) 
-\Big(\epsilon_h({f}-{u}_{tt},\hat{\theta})-\bG_t(\rv u)(\hat{\theta})\Big),
\end{eqnarray}
and hence, 
\begin{eqnarray*}
I_1=\left(\epsilon_h(\hat{f}-\hat{u}_{tt},\hat{\theta})-\bG(\rv u)(\hat{\theta})\right) 
- \int_0^t \left(\epsilon_h({f}-{u}_{tt},\hat{\theta})-\bG_s(\rv u)(\hat{\theta})\right) ds. 
\end{eqnarray*}
A use of  (\ref{ne}) for $j=1$ shows   
\begin{eqnarray}
|I_1|&\leq & |\epsilon_h(\hat{f}-({u}_{t}-u_1),\hat{\theta})|+|\bG(\rv u)(\hat{\theta})| \nonumber \\
&&+ \int_0^t \left(|\epsilon_h({f}-{u}_{tt},\hat{\theta})|+|\bG_s(\rv u)(\hat{\theta})|\right)\;ds \nonumber \\
&\leq& C \left[h^2 \left( \|\hat{f}\|_1+\|u_{t}\|_1+ \|u_1\|_1\right) + \|\bG(\rv u)\|_{-1,h}\right]\|\hat\theta\|_1\nonumber \\
&& +C  \int_0^t \left( h^2(\|f\|_1+ \|u_{tt}\|_1) +\|\bG_s(\rv u)\|_{-1,h}
\right) \|\hat\theta\|_1\; ds. 
\end{eqnarray}
Notice that $I_2$ can be written as 
\begin{eqnarray} 
I_2&=&\int_{0}^{t}\epsilon_A(\hat{\theta},\theta)ds+
\int_{0}^{t} \int_{0}^{s}\epsilon_B(\tau,\tau;\hat{\theta}(\tau),\theta(s))\;d\tau \;ds\\
&&-\int_{0}^{t}\int_{0}^{s}\int_{0}^{\tau}
\epsilon_{B_{\tau'}}(\tau,\tau';\hat{\theta}(\tau'),\theta(s)) \;d\tau'\; d\tau \;ds\\
&=& I_{21}+I_{22}+I_{23}.
\end{eqnarray}
For $I_{21}$, we apply (\ref{ng}) and the inverse inequality (\ref{inv}) to
find that 
\begin{equation} \label{newkeyy2}
|I_{21}|=\int_{0}^{t}|\epsilon_A(\hat{\theta},\theta)|ds\leq C h \int_{0}^{t}\|\theta\|_1\|\hat{\theta}\|_1\leq 
CC_{inv}\int_{0}^{t}\|\theta\|_0\; \|\hat{\theta}\|_1.
\end{equation}
In  order to estimate $I_{22}$,  we integrate by parts in time so that
\begin{eqnarray*}
|I_{22}|&=&\left|\int_0^t \epsilon_B(s,s;\hat{\theta}(s),\hat{\theta}(t))\;ds-\int_0^t
\epsilon_B(s,s;\hat{\theta}(s),\hat{\theta}(s))\;ds\right|\\
&\leq & Ch\left\{\|\hat{\theta}(t)\|_1\int_0^t\|\hat{\theta}(s)\|_1\;ds+
\int_0^t\|\hat{\theta}(s)\|_1^2\;ds\right\}.
\end{eqnarray*}
Similarly for $I_{23}$, we note that
\begin{eqnarray*}
|I_{23}|&=& \left|\int_0^t\int_0^s \epsilon_{B_{\tau}}(s,\tau;\hat{\theta}(\tau),\hat{\theta}(t))\;d\tau \;ds
-\int_0^t\int_0^s \epsilon_{B_{\tau}}(s,\tau;\hat{\theta}(\tau),\hat{\theta}(s))\;d\tau \;ds\right|\\
&\leq & 
C(T)h\left\{\|\hat{\theta}(t)\|_1\int_0^t\|\hat{\theta}(s)\|_1\;ds+
\int_0^t\|\hat{\theta}(s)\|_1^2\;ds\right\}.
\end{eqnarray*}
Using stability of $\Pi_h^*$ (i.e., $\|\Pi_h^*\theta\|_0\leq C \|\theta\|_0$) 
and the Cauchy-Schwarz inequality, it follows that 
\begin{equation}
|I_3| \leq  \int_{0}^{t}|(\rho_{t},\Pi_h^*\theta)|ds \leq 
C \int_{0}^{t}\|\rho_t(s)\|_0\|\theta(s)\|_0 ds. 
\end{equation}
For $I_4$, we apply (\ref{1.6}) and $\|\Pi_h^*\hat \theta\|_0\leq
C\|\hat \theta\|_1$ 
to obtain 
\begin{equation}
|I_4|\leq  \|u_1-\Pi_hu_1\|_0\;\|\Pi_h^*\hat{\theta}\|_0 \leq C h^2 \|u_1\|_2\;\|\hat \theta\|_1.
\label{keyy2}
\end{equation}
Finally, similarly for $I_{22}$ and $I_{23}$, an integration by parts leads to
\begin{equation} \label{ns3}
|I_5|+|I_6|\leq C(T)\left\{\|\hat{\theta}(t)\|_1\int_0^t\|\hat{\theta}(s)\|_1\;ds+
\int_0^t\|\hat{\theta}(s)\|_1^2\;ds\right\}.
\end{equation}
Now, define $\cE_0^2(t)=\|\theta(t)\|^2_0+\|\hat{\theta}(t)\|^2_1$ and let $t^\ast \in [0,t]$ 
be such that
$$\cE_0(t^\ast)=\max_{0\leq s\leq t}\cE_0(t).$$
At $t=t^\ast$, substitute the estimates (\ref{keyy1})-(\ref{ns3}) in (\ref{nn1}) and use 
the equivalence of the norms $\||\cdot\||$ and $\|\cdot\|_0$ from (\ref{nd}) along with the 
coercivity property (\ref{eqn2.3})  of $A(\cdot,\cdot)$. 
Then   a standard  use of kick back arguments  yields
\begin{eqnarray*}
\cE_0(t^\ast) &\leq& C\|\theta(0)\|+
Ch^2\left[\|u_0\|_2+\|u_1\|_2 
+\int_{0}^{t^\ast}\left(\|f\|_1+\|u_{tt}\|_1\right)\;ds\right]\\
&&+C\left[\|\bG(\rv u)\|_{-1,h}+\int_0^{t^\ast}
(\|\rho_{t}\|_0+\|\bG_s(\rv u)\|_{-1,h})\;ds\right]
+C\int_{0}^{t^\ast} \cE_0(s)\;ds.
\end{eqnarray*}
Note that $\|\theta(0)\|_0\leq Ch^2 \|u_0\|_2$.
Now apply Lemmas~\ref{lm-n}, \ref{lem2} along with  the estimates in 
Lemma~\ref{lem1} to obtain
\begin{eqnarray*}
\cE_0(t^\ast) &\leq& 
Ch^2\left( \|u_0\|_3+\|u_1\|_2 
+\int_{0}^{t^\ast}\left(\|f\|_1+\|f_{t}\|_0+\|f_{tt}\|_0\right)\;ds\right)\\
&&
+C\int_{0}^{t^\ast} \cE_0(s)\;ds.
\end{eqnarray*}
Then replace $t^\ast$ by $t$ and use  Gronwall's lemma for $t\leq T $ 
to conclude that
$$
\|\theta(t)\|_0\leq C(T)h^2 \left(\|u_0\|_3+\|u_1\|_2+
+\int_{0}^{T}\left(\|f\|_1+\|f_{t}\|_0+\|f_{tt}\|_0\right)ds \right).
$$
Finally, a  use of the triangle inequality completes the  proof.
\hfill{{\rule{2.5mm}{2.5mm}}}

\begin{remark}
 Note that it is possible to choose $u_{h,t}(0)$ as the
 $L^2$- projection of $u_1$ onto $\vv$ and in that case, the term 
 $(u_t(0)-u_{h,t}(0),\Pi_h^*\theta_t)$ becomes zero.
  \end{remark}


\subsection{Maximum norm estimates}
In this subsection, a superconvergent result for $\|\theta\|_1$ is first derived and it is then used
to analyze  quasi-optimal maximum error estimates.

\begin{lemma}
\label{H3}
Assume that $f\in L^1(H^2),\;f_t\in L^1(H^1),\;f_{tt},\;f_{ttt}\in L^1(L^2),\; u_0\in H^4\cap H_0^1$
 and $u_1\in H^3\cap H_0^1$. With $u_h(0)=\rv u_0$ and $u_{h,t}(0)=\Pi_hu_1$, there exists a positive 
 constant $C=C(T),$ independent of $h,$ such that  the following holds for $t\in (0,T]$ 
$$\|\theta_t(t)\|_0+\|{\theta}(t)\|_1\leq C\;h^2 \Big(\|u_0\|_4+\|u_1\|_3+\|D_t^3f\|_{L^1(L^2)}+
\sum_{j=0}^{2} \|D_t^jf\|_{L^1(H^{2-j})}\Big). $$
\end{lemma}
\noindent Proof. We now modify the estimates of $J_1$, $J_2$, $J_3$ and $J_5$ in (\ref{1.27}) of Theorem~ \ref{H2} to obtain a superconvergence result for $\|\theta(t)\|_1$ norm. As
$u_h(0)=\rv  u_0,$ it follows that  $ A(\theta(0),\theta(0))=0.$ 
Now with  $u_{h,t}(0)=\Pi_hu_1$, we  obtain 
\begin{equation}
|J_1|\leq Ch^4\|u_1\|_2^2.
\label{1.31}
\end{equation}
To estimate $J_2$, observe that 
$$\epsilon_h(f-u_{tt},\theta_t)=\frac{d}{dt}\epsilon_h(f-u_{tt},\theta)-\epsilon_h(f_t-u_{ttt},\theta),$$
and thus,   rewrite  $J_2$ as 
$$J_2= \epsilon_h(f-u_{tt},\theta)- \int_0^t \epsilon_h(f_t-u_{ttt},\theta)ds.$$
Then, a use of  (\ref{ne}) yields 
\begin{eqnarray}
|J_2| \leq Ch^2\Big[\left(\|f\|_1+\|u_{tt}\|_1\right)\|\theta\|_1+\int_{0}^{t}\left(\|f_t\|_1+\|u_{ttt}\|_1\right)\|\theta\|_1ds\Big].
\end{eqnarray}
For $J_3$, rewrite $G$ term as 
$$G(\rv u)(\theta_t)=\frac{d}{dt}\left\{G(\rv u)(\theta) \right\}-G_t(\rv u)(\theta),$$
and hence, a use of  (\ref{nf}) shows that 
\begin{eqnarray}
|J_3| &\leq &  |G(\rv u)(\theta)|+ \int_0^t|G_s(\rv u)(\theta)|\;ds \nonumber \\&\leq&
2\left[\|G(\rv u)\|_{-1,h}+ \int_{0}^{t}\|G_s(\rv u)\|_{-1,h}\;ds\right]\|\theta\|_1.
\end{eqnarray}
For $J_5$, apply  (\ref{1Rh}) to obtain 
\begin{equation}
|J_5|\leq 2 \int_{0}^{t}\|\rho_{tt}\|_0\|\theta_t\|_0 ds \leq C(T)h^2\int_{0}^{t}
\left(\|u_{tt}\|_2+|u_{t}\|_2 +\|u\|_2\right)\,\|\theta_t\|_0\,ds.
\label{1.32}
\end{equation}
Substituting the estimates (\ref{1.31})-(\ref{1.32}) in (\ref{1.27}), and apply standard 
kick back arguments  to arrive at 
\begin{eqnarray} \label{E-1}
\cE_1(t) &\leq&Ch^2\Big[\|u_1\|_2+\|f\|_1+\|u_{tt}\|_1+|u_{t}\|_2+\|u\|_2 \nonumber \\
&& +\int_{0}^{t}\left(\|f_t\|_1+ \|u\|_2+\|u_t\|_2+\|u_{tt}\|_2+\|u_{ttt}\|_1\right)ds\Big] \nonumber \\
 &&+C\int_{0}^{t} \cE_1(s)\;ds.\nonumber
\end{eqnarray}
An application of the integral identity (\ref{phi-t})  shows 
$$\|f\|_1\leq \|f(0)\|_1+\int_{0}^{t}\|f_t\|_1ds.$$
Then using the estimates in  Lemma \ref{lem1} we arrive at
\begin{eqnarray*}
\cE_1(t) &\leq&~Ch^2\displaystyle \Big(\|u_1\|_3+\|u_0\|_4+\|f(0)\|_1  \nonumber\\
&&+\int_{0}^{T}\left(\|f\|_2+\|f_t\|_1+ \|f_{tt}\|_0+\|f_{ttt}\|_0\right)ds\Big)\nonumber \\
 &&+C\int_{0}^{t} \cE_1(s)\;ds.
\end{eqnarray*}
Since $W^{1,1}([0,T];H^1)$ is continuously imbedded in $C^0([0,T];H^1)$, that is $\|f(0)\|_1\leq 
C\|f\|_{W^{1,1}(H^1)}$, a use of Gronwall's lemma  completes  the rest of the proof.
%
%
%
%
\hfill{{\rule{2.5mm}{2.5mm}}}
\begin{remark}
As a result of Lemma~$\ref{H3}$, we obtain a super-convergence estimate for $\theta$ in $H^1$-norm.
\end{remark}
For  $\|\theta\|_\infty,$  a  use of Sobolev inequality 
\begin{eqnarray}
\|\chi\|_\infty \leq C\left(\log\frac{1}{h}\right)^{1/2}\|\nabla\chi\|~~~~~~\forall \chi\in U_h
\label{NEWA1}
\end{eqnarray}
with  Lemma \ref{H3} yields 
\begin{eqnarray}
\|{\theta}\|_\infty &\leq & C(T)h^2 \left(\log\frac{1}{h}\right)^{1/2} \Big[\|u_0\|_4+\|u_1\|_3 \nonumber\\
&&+\int_{0}^{T}\left(\|f\|_2+\|f_t\|_1+\|f_{tt}\|_0+\|f_{ttt}\|_0\right)ds\Big].
\label{h3}
\end{eqnarray}
Below, we discuss the maximum norm estimate in form of a theorem.
\begin{theorem}\label{H4n}
Let $u$ and $u_h$  be the solutions of $(\ref{a})$ and $(\ref{c})$  respectively. Further, let the 
assumptions of  Lemma~$\ref{H3}$ hold. Then, 
\begin{eqnarray*}
\|u(t)-u_h(t)\|_\infty &\leq &C \;h^2 \left(\log\frac{1}{h}\right)\; \Big(\|u_0\|_4+\|u_1\|_3+
\|D_t^3 f\|_{L^1(L^2)}\nonumber \\ 
&&+ \sum_{j=0}^{2} \|D_t^j f\|_{L^1(H^{2-j})}\Big), 
\end{eqnarray*}
where $C=C(T)$ is a positive constant independent of $h$.
\end{theorem}
\noindent Proof.  By the triangle inequality 
$$\|u(t)-u_h(t)\|_\infty\leq \|{\theta}\|_\infty+\|{\rho}\|_\infty.$$
Now, combine the estimates obtained in (\ref{h3}) and 
in (\ref{R1h}) with Lemma~\ref{lem1} to obtain  the required result.
\hfill{{\rule{2.5mm}{2.5mm}}}



\section{Error Estimates for a Completely Discrete Scheme}
\se
In this section, we  introduce further notations and formulate 
a completely discrete scheme by applying an explicit finite difference method to 
discretize the time variable of the semidiscrete system (\ref{c}).
Then, we discuss optimal error estimates.

Let $\dt$ $(0<\dt<1)$ be the time step, $k=T/N$ for some positive integer $N$,  and $t_n=n\dt$.
For any function $\phi$ of time, let $\phi^n$ denote $\phi(t_n)$.
We shall use this notation for functions
defined for continuous in time as well as those defined for discrete in time.
Set $\phi^{n+1/2}=(\phi^{n+1}+\phi^n)/2,$
and define the following notations for the difference quotients:
$$
\pc \phi^n=\frac{\phi^{n+1}-\phi^{n-1}}{2k},
\quad \ph \phi^{n+1/2}=\frac{\phi^{n+1}-\phi^n}{k},
\quad \ps \phi^n=\frac{\phi^{n+1}-2\phi^n+\phi^{n-1}}{k^2}.
$$
Note that
$$
 \pc \phi^n=\frac{\ph \phi^{n+1/2}+\ph \phi^{n-1/2}}{2},
\qquad \ps \phi^n=\frac{\ph \phi^{n+1/2}-\ph \phi^{n-1/2}}{k}.
$$
Then, the discrete-in-time scheme of
(\ref{c}) is to seek $U^{n}\in U_h$ such that  for $\chi\in U_h$
\begin{eqnarray}
&&\frac{2}{\dt}(\ph U^{{1/2}},\Pi_h^*\chi)+A_h(U^0,\Pi_h^*\chi)=(f^0+\frac{2}{\dt}u_1,\Pi_h^*\chi), 
\quad \forall \chi\in U_h,\label{7-1}\\
&&(\ps U^{n},\Pi_h^*\chi)+A_h(U^n,\Pi_h^*\chi)+
k\sum_{j=0}^{n-1}B_h(t_n,t_{j+1/2};U^{j+1/2},\Pi_h^*\chi)=(f^n,\Pi_h^*\chi),\label{7-1c}
\end{eqnarray}
$n\geq 1$, with a given initial data $U^{0}$ in $U_h$. 
This choice of time  discretization leads to a second order accuracy in $\dt$.
The integral term in (\ref{c}) is  computed by using
the second order quadrature formula
$$
\sigma^n(g)=k\sum_{j=0}^{n-1}g(t_{j+1/2})\approx\int_0^{t_n}g(s)\,ds,\quad
\mbox{with} \quad t_{j+1/2}=(j+1/2)\dt.
$$
We shall use a shorthand notation $\sigma^n(B_h^n(U,\Pi_h^*\chi))$ for 
$k\sum_{j=0}^{n-1}B_h(t_n,t_{j+1/2};U^{j+1/2},\Pi_h^*\chi)$.
The quadrature error $q^n(g)$ is defined by
$$
q^n(g)=\sigma^{n}(g)-\int_0^{t_n}g(s)\,ds=
\sum_{j=0}^{n-1}\left( kg^{j+{1/2}}-\int_{t_j}^{t_{j+1}}g(s)\,ds \right).
$$
Similarly,  for ${\phi}\in U_h$, we define a linear functional 
$q_B^{n}({\phi})$  representing the error in 
 the quadrature formula by
$$
q_B^{n}({\phi})({\chi})=\sigma^{n}\left(B^{n}({\phi},
{\chi})\right)-\int_0^{t_n}B(t_n,s;{\phi}(s),{\chi})\,ds.
$$
Notice that $q_B^{0}({\phi})=0$. 

For our future use, we state without proof the following lemma. For a proof, see,
\cite{PTW}.
\begin{lemma}\label{q}
There exists a positive constant $C,$ independent of $k$ and $h,$ such that 
the following estimate holds:
$$
k\sum_{n=0}^{m}||\ph q_B^{n+1/2}({\phi})||_{-1,h}\leq Ck^2
\int_0^{t_{m+1}}(||{\phi}||_{1}+||{\phi}_t||_{1}+
||{\phi}_{tt}||_{1})\,ds.
$$
\end{lemma}
Now, define $e^n:=u^n-U^n$. We split $e^n=\rho^n+\xi^n$ with $\rho^n=u^n-V_hu^n$ and
$\xi^n=V_hu^n-U^n$. 
From  (\ref{7-1})-(\ref{7-1c})  and (\ref{c}), we derive equations in $e^n$ as follows
\begin{eqnarray}
&&\frac{2}{\dt}(\ph e^{{1/2}},\Pi_h^*\chi)+A_h(e^0,\Pi_h^*\chi)=(2r^0,\Pi_h^*\chi),\label{8-0-e}\\
&&(\ps e^n, \Pi_h^*\chi)+A_h(e^n,\Pi_h^*\chi)+
\sigma^n\left(B_h^n(e,\Pi_h^*\chi\right)
=(r^n,\Pi_h^*\chi)+q_{B_h}^n(u)(\Pi_h^*\chi),\label{8-3-e}
\hspace{-2cm}
\end{eqnarray}
$n\geq 1$, for all $\chi\in U_h$, where 
$\displaystyle r^0= \frac{1}{\dt}\left(\ph u^{1/2}-u_1\right)-
\frac{1}{2}u_{tt}^0=\frac{1}{2k^2}
\int_{0}^{k}(t-k)^2\frac{\partial^3 u}{\partial t^3}(t)\,dt,$ and 
\begin{equation}\label{qqw1}
r^n=\ps u^n-u_{tt}^{0}=-\frac{1}{6k^2}\\
\int_{-k}^{k}(|t|-k)^3\frac{\partial^4 u}{\partial t^4}
(t^n+t)\,dt,\;\;n\geq 1.
\end{equation}
Since estimates for $\rho$ are known from Lemma~\ref{lem2}, 
it is sufficient to estimate $\xi$. From (\ref{8-0-e})-(\ref{8-3-e}), 
we obtain the following  equations in $\xi^n$:
\begin{eqnarray}
\frac{2}{\dt}(\ph \xi^{1/2},\Pi_h^*\chi)&+&A(\xi^0,\chi)=
-\frac{2}{\dt}(\ph \rho^{{1/2}},\Pi_h^*\chi)+(2r^0,\Pi_h^*\chi) \nonumber\\
&+&\epsilon_h(f^0-u_{tt}^0,\chi)-
\epsilon_A(V_hu (0),\chi),\label{8-0}\\
(\ps \xi^n, \Pi_h^*\chi)&+&A(\xi^n,\chi)=
(r^n,\Pi_h^*\chi)-(\ps\rho^n, \Pi_h^*\chi)+H^n(\xi)(\chi) \nonumber \\
&-&\sigma^n\left(B^n(\xi,\chi)\right)- H^n(V_hu)(\chi)+\epsilon_h(f^n-u_{tt}^n,\chi)+
q^n_B(V_hu)(\chi),  \label{8-3} 
\end{eqnarray}
where 
$$H^n(\xi)(\chi)=\epsilon_A(\xi^n,\chi)+k\sum_{j=0}^{n-1}
\epsilon_B(t_n,t_{j+1/2};\xi^{j+1/2},\chi).$$

Below, we shall obtain $l^\infty(H^1)$-estimate for $\xi^{n+1/2}$.
\begin{lemma}\label{H5} 
Assume that $f\in L^1(H^2),\;f_t\in L^1(H^1),\;f_{tt},\;f_{ttt}\in L^1(L^2),\; u_0\in H^4\cap H_0^1$
 and $u_1\in H^3\cap H_0^1$. Further, assume that the CFL condition
\begin{equation}\label{CFL}
 \frac{k^2}{h^2}\leq \frac{4c_{eq}}{\Lambda C_{inv}}
\end{equation}
is satisfied, where $\Lambda>0$ is the constant given in $(\ref{eqn2.3})$, $C_{inv}$ appears in the 
inverse inequality $(\ref{inv})$ and  $c_{eq}$ is stated in the equivalence of norms as in 
$(\ref{nd})$. Then, with $u_h(0)=\rv u_0$ and $u_{h,t}(0)=\Pi_hu_1,$
there exists a positive constant $C=C(T),$ independent of $h$ and $k$, 
such that the following estimate 
\begin{eqnarray}\label{eu2-d-1}
\|\ph\xi^{m+1/2}\|_0+\|\xi^{m+1/2}\|_1
&\leq & 
C(T)(k^2+h^2) \Big(\|u_0\|_4+\|u_1\|_3\nonumber\\
&&+\|D_t^3f\|_{L^1(L^2)}+
\sum_{j=0}^{2} \|D_t^jf\|_{L^1(H^{2-j})}\Big),
\end{eqnarray}
holds for $m=0,1,\cdots,N-1$. 
\end{lemma}
\noindent Proof. Choose $\chi= \pc\xi^n$ in (\ref{8-3}) and obtain
\begin{eqnarray}\label{8-main}
\frac{1}{2} \pb \Big(|||\ph \xi^{n+{1/2}}|||_0^2 &+& 
A(\xi^{n+1},\xi^{n})\Big) = (r^n-\ps\rho^n, \delta_t
\xi^{n})+H^n(\xi)(\delta_t \xi^n)\nonumber\\
&& -\sigma^n\left(B^n(\xi,\delta_t \xi^n)\right)- H^n(V_hu)(\delta_t \xi^n)\nonumber\\
&&+\epsilon_h(f^n-u_{tt}^n,\delta_t \xi^n)+q^n_B(V_hu)(\delta_t \xi^n)\\
&=& I_1^n+ I_2^n + I_3^n+I_4^n+ I_5^n + I_6^n, \nonumber
\end{eqnarray}
where $\pb$ denotes backward differencing. 
Next multiply (\ref{8-main}) by $2 k$ and sum the resulting one from $n=2$ to $m$ to arrive at 
\begin{eqnarray}\label{8-main-1}
\frac{1}{2} \Big( |||\ph \xi^{m+1/2}|||_0^2+A(\xi^{m+1},\xi^{m})\Big)
&\leq& \frac{1}{2} \Big(|||\ph \xi^{3/2}|||_0^2+ A(\xi^{2},\xi^{1}) \Big)\nonumber\\
&&+\dt\left|\sum_{n=2}^m(I_1^n+ I_2^n + I_3^n+I_4^n+ I_5^n + I_6^n)\right|.
\end{eqnarray}
Now define
$$|||\xi^{n+1/2}|||_1^2=||\ph \xi^{n+1/2}||_0^2+||\xi^{n+1/2}||_1^2,$$
and let for some $m^{\star}$ with $ 0 \le m^{\star}\leq m,$
$$|||\xi^{m^{\star}+1/2}|||_1 =\max_{0\leq n\leq m}|||\xi^{n+1/2}|||_1.$$
To estimate the sum in $I_1^n$, an application of  the Cauchy-Schwarz 
inequality yields
\begin{eqnarray*}
\dt\left|\sum_{n=2}^m I_1^n \right|
&\leq & Ck\sum_{n=2}^m\left(\|\partial^2_t \rho^n\|_0+\|r^n\|_0\right)\, 
\left(\|\ph \xi^{n+1/2}\|_0+\|\ph \xi^{n-1/2}\|_0 \right)\\
&\leq & 2C\dt \sum_{n=2}^{m}\left(\|\partial^2_t \rho^n\|_0+\|r^n\|_0\right)
|||\xi^{m^{\star}+1/2}|||_1.
\end{eqnarray*} 
For the second sum on the right hand side of (\ref{8-main-1}), 
we use the fact that
\begin{eqnarray}\label{eq:formula}
\psi^{n}\pc\xi^{n}&=&\pb(\psi^{n}\xi^{n+1/2})-\ph\psi^{n+1/2}\xi^{n-1/2}
\end{eqnarray}
and conclude 
$$
k\sum_{n=2}^m\epsilon_A(\xi^n,\pc\xi^{n})= \epsilon_A(\xi^m,\xi^{m+1/2})-
\epsilon_A(\xi^1,\xi^{1+1/2})-k\sum_{n=2}^m\epsilon_A(\ph\xi^{n+1/2},\xi^{n-1/2}).
$$
Using (\ref{ng}) and the inverse inequality (\ref{inv}),  we obtain
\begin{eqnarray*}
\left|k\sum_{n=2}^m\epsilon_A(\xi^n,\pc\xi^{n})\right|
&\leq&
Ch\left\{\|\xi^m\|_1\|\xi^{m+1/2}\|_1+\|\xi^1\|_1\|\xi^{1+1/2}\|_1\right\}\\
&& +Chk\sum_{n=2}^m\|\ph\xi^{n+1/2}\|_1 \|\xi^{n-1/2}\|_1\\
&\leq& C\left\{\|\xi^m\|_0+\|\xi^1\|_0+k\sum_{n=2}^m\|\ph\xi^{n+1/2}\|_0\right\}
|||\xi^{m^{\star}+1/2}|||_1.
\end{eqnarray*}
Since $\xi^0=0$, $\xi^m=k\sum_{n=0}^{m-1} \ph \xi^{n+1/2}$, 
and  it follows that
$$
\left|k\sum_{n=2}^m\epsilon_A(\xi^n,\pc\xi^{n})\right|\leq 
Ck\left(\sum_{n=0}^{m-1} \|\ph \xi^{n+1/2}\|_0\right)
|||\xi^{m^{\star}+1/2}|||_1.
$$
Similarly, we obtain
\begin{eqnarray*}
\left|k^2\sum_{n=2}^m\sum_{j-0}^{n-1}\epsilon_B(t_n,t_{j+1/2};\xi^{j+1/2},\pc\xi^{n})\right|
&\leq&k^2\sum_{n=2}^m\left(C\sum_{j=0}^{n-1}\|\xi^{j+1/2}\|_1\right)
|||\xi^{m^{\star}+1/2}|||_1\\
&\leq&CT k\left(\sum_{j=0}^{m-1}\|\xi^{j+1/2}\|_1\right)
|||\xi^{m^{\star}+1/2}|||_1,
\end{eqnarray*} 
and hence, 
$$
\dt\left|\sum_{n=2}^m I_2^n \right|\leq
C(T)k\left(\sum_{n=0}^{m-1} |||\xi^{n+1/2}|||_1\right)
|||\xi^{m^{\star}+1/2}|||_1.
$$
To estimate the sum in $I_3^n$, we again use (\ref{eq:formula}) and 
rewrite the sum as:
\begin{eqnarray*}
k\sum_{n=2}^mI_3^n
&=&\sigma^m\left(B^m(\xi,\xi^{m+1/2})\right)-\sigma^1\left(B^1(\xi,\xi^{1+1/2})\right)\\
&&-k^2\sum_{n=2}^m\sum_{j=0}^{n-1}(\bar{\partial}_{t,1}B)
(t_n,t_{j+1/2};\xi^{j+1/2},\xi^{n-1/2})\\
&&+k\sum_{n=2}^mB(t_{n-1},t_{n-1/2};\xi^{n-1/2},\xi^{n-1/2}),
\end{eqnarray*}
where $\bar{\partial}_{t,1}B$ denotes the difference quotient of $B$ with respect to
its first argument. Since, 
$|\bar{\partial}_{t,1}B|\leq C||B_t||_\infty< \infty$, it follows that
$$
\left|k\sum_{n=2}^mI_3^n\right|\leq
C(T)k\left(\sum_{j=0}^{m-1}||\xi^{j+1/2}||_1\right)
|||\xi^{m^{\star}+1/2}|||_1. 
$$
For the sum involving $I_4^n$, we note that
\begin{eqnarray*}
\left|k\sum_{n=2}^m\epsilon_A(V_hu^n,\pc\xi^{n})\right|&=& 
\left|\epsilon_A(V_hu^m,\xi^{m+1/2})-
\epsilon_A(V_hu^1,\xi^{1+1/2})-k\sum_{n=2}^m\epsilon_A(\ph V_hu^{n+1/2},\xi^{n-1/2})\right|\\
&\leq& Ch^2\left\{||u^m||_2+||u^1||_2+k\sum_{n=0}^m||\ph V_hu^{n+1/2}||_2\right\}
|||\xi^{m^{\star}+1/2}|||_1\\
&\leq& Ch^2\left\{ ||u_0||_2+||u_t||_{L^1(H^2)}\right\}
|||\xi^{m^{\star}+1/2}|||_1.
\end{eqnarray*}
Similarly, we  have
$$
\left|k^2\sum_{n=2}^m\sum_{j=0}^{n-1}\epsilon_B(t_n,t_{j+1/2};V_hu^{j+1/2},\pc\xi^{n})\right|
\leq CTh^2\left\{ ||u_0||_2+||u_t||_{L^1(H^2)}\right\} |||\xi^{m^{\star}+1/2}|||_1.
$$
In order to estimate the sum in $I_5^n$, we repeat the previous arguments 
and use (\ref{ne}) to arrive at
\begin{eqnarray*}
\left|k\sum_{n=2}^m\epsilon_h(f^n-u^n_{tt},\pc\xi^{n})\right|&=& 
\left|\epsilon_h(f^m-u^m_{tt},\xi^{m+1/2})-
\epsilon_h(f^1-u^1_{tt},\xi^{1+1/2})\right.\\
&&\left.-k\sum_{n=2}^m\epsilon_h\left(\ph (f^{n+1/2}-u^{n+1/2}_{tt}),\xi^{n-1/2}\right)\right|\\
&\leq& Ch^2\left\{ \|f^0-u_{tt}^0\|_1+ \|f_t-u_{ttt}\|_{L^1(H^1)}\right\}
|||\xi^{m^{\star}+1/2}|||_1.
\end{eqnarray*}
For the last sum, we rewrite it as 
$$
k\sum_{n=2}^mI_6^n=q^m_B(V_hu)(\xi^{m+1/2})-q^1_B(V_hu)(\xi^{1+1/2})
-k\sum_{n=2}^m(\ph q^{n+1/2}_B(V_hu))(\xi^{n-1/2}).
$$
Since $q^0_B(V_hu)=0$, $q^m_B(V_hu)=k\sum_{n=0}^m\ph q^{n+1/2}_B(V_hu)$, 
we obtain
$$
k\left|\sum_{n=2}^mI_6^n\right|\leq Ck\left\{\sum_{n=0}^m
||\ph q^{n+1/2}_B(V_hu)||_{-1,h}\right\}|||\xi^{m^{\star}+1/2}|||_1.
$$
Combining all the previous estimates, we conclude that
\begin{eqnarray}\label{8-main-2}
|||\ph \xi^{m+1/2}|||_0^2+A(\xi^{m+1},\xi^{m})
&\leq&|||\ph \xi^{3/2}|||_0^2+ A(\xi^{2},\xi^{1})+
Ck\left\{ \sum_{n=2}^{m}(\|\partial^2_t \rho^n\|_0+\|r^n\|_0)\right.\nonumber\\
&+& \left.\sum_{n=0}^m||\ph q^{n+1/2}_B(V_hu)||_{-1,h}+
\sum_{j=0}^{m-1}|||\xi^{j+1/2}|||_1\right\}|||\xi^{m^{\star}+1/2}|||_1\nonumber\\
&+& h^2C(T,f,u)|||\xi^{m^{\star}+1/2}|||_1,
\end{eqnarray}
where
$$
C(T,f,u)=||u_0||_2+||u_t||_{L^1(H^2)}+
||u_{tt}(0)||_1+||u_{ttt}||_{L^1(H^1)}
+||f^0||_1+||f_t||_{L^1(H^1)}.
$$
In order to estimate the first two terms on the right hand side of
(\ref{8-main-2}), we choose $\chi=\ph\xi^{3/2}$ in (\ref{8-3}) for $n=1$ and
obtain
\begin{eqnarray*}
|||\ph \xi^{3/2}|||_0^2+A(\xi^{2},\xi^{1})
&\leq&|||\ph \xi^{1/2}|||_0^2+h^2\Big(\|u^1\|_2+\|u_0\|_2+k\|\ph u^{1/2}\|_2\Big)\\
&+&h^2\Big( \|f^0-u_{tt}^0\|_2+\|f_t-u_{ttt}\|_{L^1(0,k;H^1)} \Big) +\|\ph q^{1/2}_B\|_{-1,h}.
\end{eqnarray*}
Next, we choose  $\chi=\ph\xi^{1/2}$ in (\ref{8-0}) to find that
$$
|||\ph \xi^{1/2}|||_0\leq C\left\{\|\ph \rho^{1/2}\|_0+k\|r^0\|_0+h^2\|f^0-u_{tt}^0\|_2+
h^2\|u_0\|_2\right\}.
$$ 
A use of  these estimates in (\ref{8-main-2}) results in 
\begin{eqnarray}\label{8-main-3}
|||\ph \xi^{m+1/2}|||_0^2&+&A(\xi^{m+1},\xi^{m})
\leq
C\left\{\|\ph \rho^{1/2}\|_0+k\sum_{n=1}^{m}\|\partial^2_t \rho^n\|_0+
k\sum_{n=0}^{m}\|r^n\|_0\right.\nonumber\\
&&\left.+k\sum_{n=0}^m||\ph q^{n+1/2}_B(V_hu)||_{-1,h}+
k\sum_{j=0}^{m-1}||\xi^{j+1/2}||_1\right\}|||\xi^{m^{\star}+1/2}|||_1\nonumber\\
&&+h^2C(T,f,u)|||\xi^{m^{\star}+1/2}|||_1.
\end{eqnarray}
Note that
$$A(\xi^{m+1},\xi^{m})=A(\xi^{m+1/2},\xi^{m+1/2})-
\frac{k^2}{4}A(\ph\xi^{m+1/2},\ph\xi^{m+1/2}).$$
Hence,
$$
|||\ph \xi^{m+1/2}|||_0^2+A(\xi^{m+1},\xi^{m})\geq
c_{eq} \;\|\ph \xi^{m+1/2}\|_0^2 + \alpha\|\xi^{m+1/2}\|_1^2-
\frac{k^2}{4}A(\ph\xi^{m+1/2},\ph\xi^{m+1/2}).
$$
Since the CFL condition (\ref{CFL}) holds, choose $k$ so that $C_\ast=\left(c_{eq}-\Lambda C_{inv} \frac{k^2}{4h^2}\right)>0$, 
where the constants $\Lambda$, $c_{eq}$ and $C_{inv}$ appear in (\ref{eqn2.3}), 
(\ref{nd}) and (\ref{inv}), respectively.  Then
$$
|||\ph \xi^{m+1/2}|||_0^2+A(\xi^{m+1},\xi^{m})\geq \min\{C_\ast,\alpha\}
|||\xi^{m+1/2}|||_1.
$$
Altogether, it now results in
\begin{eqnarray}\label{8-5}
|||\xi^{m+1/2}|||_1\leq|||\xi^{m^{\star}+1/2}|||_1
&\leq&
C\left\{\|\ph \rho^{1/2}\|_0+k\sum_{n=1}^{m}\|\partial^2_t\rho^n\|_0+
k\sum_{n=0}^{m}\|r^n\|_0\right.\nonumber\\
&&\left.+k\sum_{n=0}^m||\ph q^{n+1/2}_B(V_hu)||_{-1,h}+
k\sum_{j=0}^{m-1}|||\xi^{j+1/2}|||_1\right\}\nonumber\\
&&+h^2C(T,f,u).
\end{eqnarray}
To estimate the first two terms on the right hand side of (\ref{8-5}), it is observed that
\begin{eqnarray}\label{8-5a}
||\ph\rho^{1/2}||_0\leq\frac{1}{k}\int_0^{\dt}||\rho_t(s)||_0\,ds,
\end{eqnarray}
and a use of Taylor series expansion yields
\begin{eqnarray}\label{8-5b-n}
k\sum_{n=1}^{m}||\ps\rho^n||_0 &\leq&\frac{1}{k}\sum_{n=1}^{m}\left\{
\int_{t_n}^{t_{n+1}}(t_{n+1}-s)||\rho_{tt}(s)||_0\,ds+
\int_{t_{n-1}}^{t_n}(s-t_{n-1})||\rho_{tt}(s)||_0\,ds
\right\}\nonumber\\
&\leq& 2\int_{0}^{t_{m+1}}||\rho_{tt}(s)||_0 \,ds.
\end{eqnarray}
Further, from (\ref{qqw1}) it follows that
$$||r^n||_0\leq Ck\int_{t_{n-1}}^{t_{n+1}}\|D_t^4u(s)\|_0\,ds,\quad n\geq 1,$$
and 
$$||r^0||_0\leq Ck||u_{ttt}||_{L^\infty(0,\dt/2;L^2(\Omega))}
\leq Ck\int_{0}^{t_{m+1}} \|D_t^3u(s)\|_0 \,ds.$$
Thus, we arrive at
\begin{equation}\label{s2}
k\sum_{n=0}^{m}||r^n||_0\leq Ck^2\int_{0}^{t_{m+1}}
\left( \|D_t^3u(s)\|_0+\|D_t^4u(s)\|_0\right)\,ds.
\end{equation}
Finally, a use of Lemma~\ref{q} and the triangle
inequality yields
$$
k\sum_{n=1}^{m}||\ph q^{n+1/2}_B(V_hu)||_{-1,h}\leq Ck^2
\sum_{j=0}^2\int_{0}^{t_{m+1}} \left(\|D_t^j u(s)\|_1+\|D_t^j \rho(s)\|_1\right)\,ds.
$$
Substitute now (\ref{8-5a})-(\ref{s2}) in (\ref{8-5}) and use the estimates in 
Lemmas~\ref{lem2} and \ref{lem1}. Then, an application  of
the discrete Gronwall's lemma completes the 
rest of the proof. 
\hfill{{\rule{2.5mm}{2.5mm}}}

By Sobolev inequality, it follows that 
\begin{equation}\label{s4}
\|\xi^{n+1/2}\|_\infty \leq C\left(\log \frac{1}{h}\right)^{1/2} \|\xi^{n+1/2}\|_1.
\end{equation}
Using Lemma~$\ref{H5}$, the triangle inequality and 
the estimates (\ref{s4}) and (\ref{R1h}), we obtain the result of th following theorem.
\begin{theorem}\label{H6}
Let the  assumptions of  Lemma~$\ref{H5}$ hold. Then, 
\begin{eqnarray}\label{eu2-d-2}
\|u(t_{m+1/2})-U^{m+1/2}\|_\infty
&\leq & 
C(T)\left(\log \frac{1}{h}\right)(k^2+h^2) \Big(\|u_0\|_4+\|u_1\|_3\nonumber\\
&&+\|D_t^3f\|_{L^1(L^2)}+
\sum_{j=0}^{2} \|D_t^jf\|_{L^1(H^{2-j})}\Big)
\end{eqnarray}
for $m=0,1,\cdots,N-1$.
\end{theorem}

\section{ FVEM with Quadrature}
\se
In this section, we discuss the effect of numerical quadrature on FVEM, when the $L^2$ 
inner product $(\cdot,\cdot)$ and the bilinear forms $A_h(\cdot,\cdot)$ and 
$B_h(t,s;\cdot,\cdot)$ appearing in (\ref{c}) are approximated  by simple quadrature formulae.

For a continuous function $\phi$ on a triangle $K$, consider the quadrature formula  
\begin{eqnarray}
\mathcal {Q}_{K,h}(\phi)=\frac{1}{3}|K|\sum_{l=1}^3 \phi(P_l) \approx \int_K\phi(x) dx\;\;~~~~\forall K \in \cT_h,
\label{Q1}
\end{eqnarray}
where $P_l,~1\leq l\leq 3$ denote the vertices of the triangle $K$
and $|K|$ denotes the area of the triangle $K$. Now the quadrature formula given by
(\ref{Q1}) is exact for $\phi \in P_1(K) ~\forall K\in \cT_h$.
Using (\ref{Q1}), we replace the $L^2$ inner product by the following discrete $L^2$ inner
product:
\begin{eqnarray}
(\chi,\Pi_h^*\psi)_h&=&\sum_{K\in \cT_h}\mathcal
{Q}_{K,h}(\chi\Pi_h^*\psi) \nonumber \\
&=&\sum_{P_i\in N_h^0} \chi(P_i)\psi(P_i)|S_{K_{P_i}}^*|~~~~~~\forall \chi,~\psi \in U_h.
\label{QQ12}
\end{eqnarray}
This is known as lumping of mass in the literature.
Observe that $\|\chi\|_h^2=(\chi,\chi)_h~~\forall \chi \in U_h$ is a
norm on $U_h,$ which is equivalent to the $L^2$ norm, i.e.,  there exist
 positive constants $C_5$ and $C_{6}$,  independent of $h$,  such that 
 \begin{equation}
 C_5\|\chi\|_0\leq\|\chi\|_h\leq C_{6}\|\chi\|_0.
 \label{ap7}
 \end{equation}
Define quadrature error  by
$$\bar{\epsilon}_h(\chi,\psi)=(\chi,\Pi_h^*\psi)-
(\chi,\Pi_h^*\psi)_h.$$ 
Since the quadrature formula involves only the
values of the functions at the interior nodes  and 
$\Pi_h^*u_h(P_i)=u_h(P_i)\;\forall P_i\in N_h^0 \;\mbox{and}\; u_h\in U_h$, 
it follows that 
\begin{eqnarray}
 (\chi,\psi )_h=(\chi,\Pi_h^*\psi )_h\; \;\;\;\forall \chi,\;\psi \in U_h.
 \label{Q10}
\end{eqnarray}
Below, we state the  estimates related to quadrature error, whose proof can be 
found in  \cite{KNP-2008}.
\begin{lemma}
\label{L21}
For  $\chi,~ \psi \in U_h$, there is a positive constant $C$, 
independent of $h$,  such that  the following estimate holds:
\begin{eqnarray}
|\bar{\epsilon}_h(\chi,\psi)|\leq Ch^2\|\chi\|_1\|\psi\|_1.
\label{apnn1}
\end{eqnarray}
Further, for $\chi \in H^2$ and $\psi \in U_h$,  there holds: 
\begin{eqnarray}
|\bar{\epsilon}_h(\chi,\psi)|\leq Ch^2\|\chi\|_2\|\psi\|_1.
\label{apnn2}
\end{eqnarray}
\end{lemma}

Now define the following quadrature approximation over each element $K$ by
\begin{eqnarray}
\label{Q2}
\int_{\overline {M_lQ}\cap K}v(z)~ds \approx
\frac {\overline{M_lQ}}{2}\left(v(M_l)+v(Q)\right)= \tilde
{\mathcal{Q}}_{h,l}(v),
\end{eqnarray}
 \begin{figure}
    \begin{center}
    \includegraphics*[width=8.0cm,height=6.0cm]{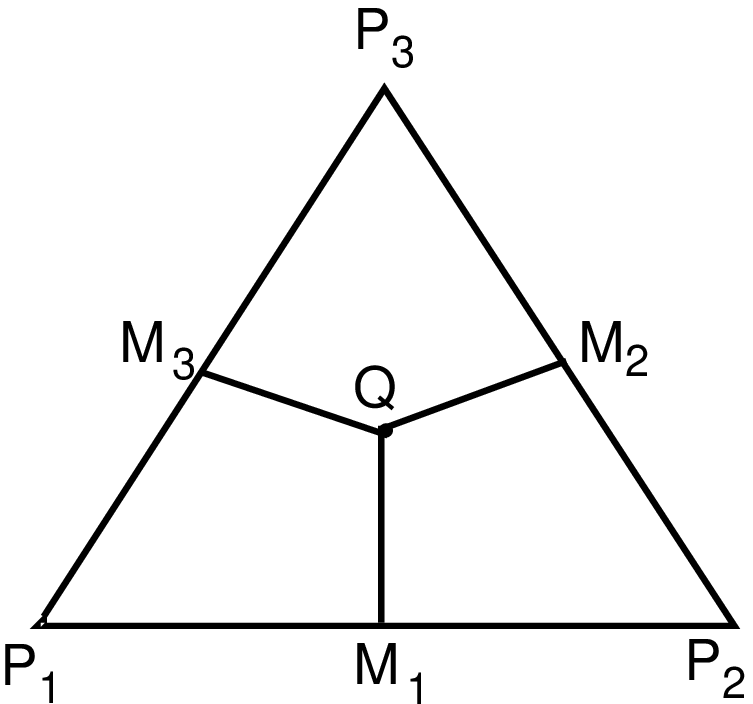}
       \caption{}
       \label{fig:mesh-1}
    \end{center}
\end{figure}
where $M_l$ is the midpoint of $P_lP_{l+1}$ and $Q$ is the barycenter
of the triangle $\triangle P_lP_{l+1}P_{l+2}$,  (see  FIGURE~\ref{fig:mesh-1}  for $l=1$).
Associated with (\ref{Q2}), we now introduce the quadrature error as
$$ \mathcal E_{\overline {M_lQ}\cap K}(v)=\int_{\overline {M_lQ}\cap K}v(s)ds- \tilde {\mathcal { Q}}_{h,l}(v).$$
Then, we have the following estimate related to the above quadrature error. For a proof,
see, Cai \cite[pp 732]{10}.
\begin{lemma}
\label{AP1}
Let $v\in W^{2,\infty}({\overline {M_lQ}\cap K}).$  Then, there is a
positive constant $C,$ independent of $h_K$,  such that 
\begin{eqnarray}
|\mathcal E_{\overline {M_lQ}\cap K}(v)|\leq Ch_K^3\|v\|_{2,\infty,\overline {M_lQ}\cap K},
\label{ni}
\end{eqnarray}
where $h_K$ is the diam($K$).
\end{lemma}
\noindent 

Now to replace the integral in the definition of $A_h(\cdot,\cdot)$,  we
observe that
\begin{eqnarray*}
A_h(u_h,\Pi^*_hv_h)&=&-\sum_{P_l\in N_h}v_i\int_{\partial K_{P_l}^* } A\nabla
u_h.{\bf n}~ds\;\;\;\Big(v_i=v_h(P_i)\Big)\\
&=& \sum_{K}I_K (u_h, \Pi_h^*v_h),
\end{eqnarray*}
where 
\begin{eqnarray*}
I_K (u_h, \Pi_h^*v_h)&=&-\sum_{{P_l}(1\leq l\leq 3)}  v_l \int_{\partial K_{P_l}^*\cap K} A\nabla
u_h.{\bf n_l} ds\\
&=&\sum_{{P_l}(1\leq l\leq 3)}(v_{l+1}-v_{l})
\int_{\overline{M_lQ}\cap K} A\nabla
u_h.{\bf n_l} ds,
\end{eqnarray*}
$v_4=v_1$ and ${\bf n_l}$ is the outward unit
normal vector to $\overline{M_lQ}$. Since $\nabla u_h.{\bf n_l}$ is constant
on each element $K$, we  define the quadrature rule  as
\begin{eqnarray}
\tilde {I}_K(u_h, \Pi_h^*v_h)=\sum_{{P_l}(1\leq l\leq 3)}  \mathcal E_{\overline {M_lQ}\cap K}(A)\nabla u_h.{\bf n_l}(v_{l+2}-v_{l+1}).
\label{Q123}
\end{eqnarray}
and set
$$\tAh(\chi,\Pi_h^*\psi)= \sum_{K\in
T_h}\tilde{I}_K(\chi,\Pi_h^*\psi).$$
Note that the bilinear form $A_h(\cdot,\cdot)$ in (\ref{c}) is  approximated by
$\tAh(\cdot,\cdot).$
Simlilarly, define $\tBh(\cdot,\cdot)$ as an approximation  of $B_h(\cdot,\cdot).$

With the definitions as above, define quadrature error functional for the bilinear
form $A_h(\cdot,\cdot)$ as
\begin{eqnarray}
\beA(\chi,\psi)=A_h(\chi,\Pi^*_h\psi)-\tAh(\chi,\Pi^*_h\psi)\;\;\;\; \forall \chi,~\psi \in U_h.
\label{Q4}
\end{eqnarray}
Below, we state without proof the estimate of (\ref{Q4}) whose proof can be found in  \cite{KNP-2008}.
\begin{lemma}
\label{L1}
Assume that $\cA \in  W^{2,\infty} (\Omega; \mathbb{R}^{2\times 2}).$ Then, there exists a positive 
constant $C,$ independent of $h,$ such that
$$\beA(\chi,\psi)
 \leq Ch^2\|\chi\|_1\|\psi\|_1\;\;\;\;\forall \chi,~\psi \in U_h.$$
\end{lemma}
Similar results hold for $\beB(t,s;\cdot,\cdot)$ which is defined as in (\ref{Q4}). 
For the rest of our analysis, we introduce the functionals $\pp(t)=\pp$ and 
$\hat{S}(t)=\hat{S}$  defined on $U_h$ for a given $\psi$ and $t\in (0,T]$ as
$$
\pp(\psi)(\chi)=\beA(\psi,\chi)+\int_0^t\beB(t,s;\psi(s),\chi)ds,
$$
and
$$
\hat{S}(\psi)(\chi)=\beA(\hat\psi,\chi)+\int_0^t\beB(s,s;\hat\psi (s),\chi)\,ds
-\int_0^t\int_0^s\bar{\epsilon}_{B_\tau}(s,\tau;\hat\psi (\tau),\chi)\,d\tau ds.
$$
Then using Lemma~\ref{L1}, we derive the following estimate for $S$  in a similar manner  to those 
obtained in Lemma~\ref{lm-n} 
$$
\|S(\psi)\|_{-1,h}\leq Ch^2  \left( \|\psi(t) \|_{2}+ \int_{0}^t \|\psi(s)\|_{2}\,ds\right).
$$
Similar result can be obtain for the estimate of $\hat{S}$ again following proof of Lemma~\ref{lm-n}.
 
Now the semidiscrete finite volume element method combined with quadrature is to seek
$u_h:(0,T]\longrightarrow U_h$ such that  
\begin{equation}
(u_{h,tt},v_h)_h+\tAh(u_h,v_h)+\int_0^t\tBh(t,s;u_h(s),v_h)\,ds=(f,v_h)_h~~~~~~\forall v_h 
\in \vv,\label{ap4}
\end{equation}
with  appropriate initial data $u_h(0)$ and $u_{h,t}(0)$ in $U_h$.

\subsection{ Optimal error estimates }
In this subsection, we discuss optimal estimates in $L^{\infty}(L^2)$  as well as in $L^{\infty}(H^1)$-norms and quasi-optimal estimates in $L^{\infty}(L^{\infty})$-norm.

Now replace $v_h$ by $\Pi_h^*\chi$ in (\ref{ap4}) 
and subtract the resulting equation from (\ref{newaa1}) to obtain 
\begin{eqnarray}
(u_{tt},\Pi_h^*\chi)-(u_{h,tt},\Pi_h^*\chi)_h&+&A_h(u,\Pi_h^*\chi)-\tAh(u_h,\Pi_h^* \chi)\nonumber \\
&+&\int_0^tB_h(t,s;u,\Pi_h^*\chi)\,ds-\int_0^t\tBh(t,s;u_h,\Pi_h^*\chi)\,ds \nonumber\\
&&= (f,\Pi_h^*\chi)-(f,\Pi_h^*\chi)_h \;\;\; \forall  \chi \in U_h.
\label{equation1}
\end{eqnarray}
Using the definitions of Ritz-Volterra projection $\rv u$ and $\pp,$  we arrive at an equation in $\theta$ as
\begin{eqnarray}
(\theta_{tt},\Pi_h^*\chi)_h+A(\theta, \chi)&=&
-(\rho_{tt},\Pi_h^*\chi)+\G(\rho)(\chi)+\G(\theta)(\chi)\nonumber \\
&&+\pp(\theta)(\chi)- \pp(\rv u)(\chi)
+\beh(f,\chi)-\beh((\rv u)_{tt},\chi)\nonumber\\
&&-\int_0^tB(t,s;\theta,\chi)\,ds.
\label{ap5-n}
\end{eqnarray}

Below, we establish $L^{\infty}(H^1)$ estimate.

\begin{theorem}
\label{HH1}
Let $u$ and $u_h$  be the solutions of $(\ref{a})$ and $(\ref{ap4}),$  respectively, 
and assume that $f\in L^1(H^1),~f_t,~f_{tt}\in L^1(L^2),~u_0\in H^3\cap H_0^1$ and 
$u_1\in H^2\cap H_0^1$. With $u_h(0)=\Pi_h u_0$ and $ u_{h,t}(0)=\Pi_h u_1,$  
there exists a positive constant $C=C(T),$ independent  of $h$,  such that 
 $$\|u(t)-u_h(t)\|_1 \leq Ch \left(\|u_0\|_3+\|u_1\|_2+ \|f\|_{L^1(H^1)}+
 \sum_{j=1}^{2} \|D_t^j f\|_{L^1(L^2)}\right) $$
 holds for $t\in (0,T].$
\end{theorem}

\noindent Proof. Choose $\chi=\theta_t$ in (\ref{ap5-n}) so that  
\begin{eqnarray}
 (\theta_{tt},\Pi_h^*\theta_t)_h+A(\theta,\theta_t)&=&
\G(\rho)(\theta_t) +\G(\theta)(\theta_t)-(\rho_{tt},\Pi_h^*\theta_t)+\pp(\theta)(\theta_t)\nonumber \\
&&-\pp(\rv u)(\theta_t)+\bar{\epsilon}_h(f,\theta_t)-\bar{\epsilon}_h((\rv u)_{tt},\theta_t)\nonumber \\
&&-\int_0^tB(t,s;\theta,\theta_t)ds.
\label{nk}
\end{eqnarray}
Then, use (\ref{Q10}) and the symmetric property of $A(\cdot,\cdot)$ to
obtain
\begin{eqnarray*}
\frac{1}{2} \frac{d}{dt}\left[(\theta_t,\theta_t)_h+A(\theta,\theta)\right] &=&
\G(\rho)(\theta_t) +\G(\theta)(\theta_t)-(\rho_{tt},\Pi_h^*\theta_t)+\pp(\theta)(\theta_t) \\
&&-\pp(\rv u)(\theta_t)+\bar{\epsilon}_h(f,\theta_t)-\bar{\epsilon}_h((\rv u)_{tt},\theta_t) \\
&&-\int_0^tB(t,s;\theta,\theta_t)ds.
\end{eqnarray*}
Integrate from $0$ to $t$ and use  the equivalence of the  norms in (\ref{ap7})
to find that 
\begin{eqnarray}
\frac{1}{2}\left[\|\theta_t\|_h^2+ A(\theta,\theta)\right]&=&\Big\{\frac{1}{2}\|\theta_t(0)\|_h^2+ 
\frac{1}{2}A(\theta(0),\theta(0)) + 
\int_0^t\Big[\G(\rho)(\theta_t)+\G(\theta)(\theta_t) \nonumber \\
&&-(\rho_{tt},\Pi_h^*\theta_t)-\int_0^sB(s,\tau;\theta(\tau),\theta_t)\,d\tau\Big]\;ds\Big\} +
 \int_0^t\pp(\theta)(\theta_t) \;ds\nonumber\\
&&- \int_0^t\pp(\rv u)(\theta_t) ds 
+ \int_0^t\bar{\epsilon}_h(f,\theta_t)ds-\int_0^t\bar{\epsilon}_h((\rv u)_{tt},\theta_t)\;ds \nonumber\\
&=&I+J_1+J_2+J_3+J_4.
\label{asp1}
\end{eqnarray}
Estimates for the first term $I$ have already been derived in Theorem \ref{H2}.
In order to estimate $J_1$, use   Lemma \ref{L1} and inverse inequality (\ref{inv}) to obtain 
\begin{eqnarray}
|J_{1}| &\leq& \int_0^t|\pp(\theta)(\theta_t)|\;ds \leq \int_0^t \|\pp(\theta)\|_{-1,h}\;\|(\theta_t)\|_1\;ds  \nonumber\\
&\leq& Ch^2\left[\int_0^t\|\theta\|_1\|\theta_t\|_1 ds +\int_0^t\int_0^s\|\theta_t(\tau)\|_1
\|\theta(s)\|_1\;d\tau \;ds\right]\nonumber \\
&\leq& Ch\left[\int_0^t\|\theta\|_1\|\theta_t\|_0 ds+\left(\int_0^t\|\theta\|_1 \;ds\right)
\left(\int_0^t\|\theta_t\|_0 \; ds\right)\right]. 
\label{asp2}
\end{eqnarray}
%
%
For $J_2$, we find that
\begin{equation}
 |J_{2}|\leq  \int_0^t|\pp(\rv u)(\theta_t)|\;ds\leq C h^{-1}
 \int_0^t\|\pp(\rv u)\|_{-1,h}\|{\theta_t}\|_0\;ds. 
\end{equation} 
In view of Lemma~\ref{L21}, the terms  $J_3$ and $J_4$ are bounded as 
\begin{eqnarray}
|J_3|+|J_4|&\leq& \int_0^t|\bar{\epsilon}_h(f,{\theta_t})|ds+2 \int_0^t|\bar{\epsilon}_h((\rv u)_{tt},{\theta_t})|ds \nonumber \\  &\leq&
Ch^2\int_0^t(\|f\|_2+\|\rho_{tt}\|_1+\|u_{tt}\|_1)\|{\theta_t}\|_1ds\nonumber \\ 
&\leq& Ch\int_0^t(\|f\|_2+\|\rho_{tt}\|_1+\|u_{tt}\|_1)\|{\theta_t}\|_0ds.
\label{asp3}
\end{eqnarray}
Now, substitute (\ref{asp2})-(\ref{asp3}) in (\ref{asp1}). Use the coercivity property of 
the bilinear form $A(\cdot,\cdot)$ and equivalence of  norms (\ref{ap7}). Then,  
proceed as in Theorem  \ref{H2} to  complete the rest of the proof.\hfill{{\rule{2.5mm}{2.5mm}}}

In the following theorem, we prove optimal $L^{\infty}(L^2)$-estimate.
\begin{theorem}
Under the assumptions of Theorem~$\ref{HH1}$,  there exists a positive constant $C=C(T)$, 
independent of $h$, such that 
$$\|u(t)-u_h(t)\|_0\leq Ch^2 \left(\|u_0\|_3+\|u_1\|_2+ \|f\|_{L^1(H^1)}+\sum_{j=1}^{2} 
\|D_t^jf \|_{L^1(L^2)}\right)$$
holds for all $t\in (0,T].$
\end{theorem}
 \noindent Proof.  
Integrate (\ref{ap5-n}) from $0$ to $t$ to arrive at
\begin{eqnarray}
(\theta_t,\Pi_h^*\chi)_h &+& A(\hat{\theta}, \chi) = -(\rho_t,\Pi_h^*\chi)+\bG(\rho)(\chi)+\bG(\theta)(\chi)\nonumber \\
&+&\bpp(\theta)(\chi)-\bpp(\rv u)(\chi)
+\bar{\epsilon}_h(\hat{f},\chi)-\bar{\epsilon}_h((\rv u)_t,\chi)\nonumber \\
&+& (u_{t}(0),\Pi_h^*\chi)-(u_{h,t}(0),\Pi_h^*\chi)_h-\int_0^t\int_0^sB(s,\tau;\theta(\tau),\chi)\,d\tau \;ds.
\label{ap5}
\end{eqnarray}
Choose  $\chi=\theta$ in (\ref{ap5}) and  use (\ref{Q10}) with the symmetry of the bilinear 
form $A(\cdot,\cdot)$ to obtain 
\begin{eqnarray}
\frac{1}{2}\frac{d}{dt}\left[(\theta,\theta)_h+A(\hat{\theta},\hat\theta)\right]&=& 
I(t)+\bpp(\theta)(\theta)-\bpp(\rv u)(\theta)
+\bar{\epsilon}_h(\hat{f},\theta)\nonumber \\
&&-\bar{\epsilon}_h((\rv u)_t,\theta)+(u_{t}(0),\Pi_h^*\theta)-(u_{h,t}(0),\Pi_h^*\theta)_h.
\label{bbbb1}
\end{eqnarray}
where 
$$
I(t)=-(\rho_t,\Pi_h^*\theta)+\bG(\rho)(\theta)+\bG(\theta)(\theta)-
\int_0^t\int_0^sB(s,\tau;\theta(\tau),\theta)\,d\tau ds.
$$
Integrate (\ref{bbbb1}) from $0$ to $t$  to  find that
\begin{eqnarray}
\frac{1}{2}\Big(\|\theta(t)\|_h^2+A(\hat{\theta},\hat\theta)\Big) &=& \frac{1}{2} \|\theta(0)\|_h^2+
\int_0^tI(s)\,ds+ \int_0^t\bpp(\theta)(\theta) \;ds \nonumber \\
&& - \int_0^t \bpp(\rv u)(\theta)\;ds
-\int_0^t \bar{\epsilon}_h((\rv u)_t,\theta)\; ds \nonumber \\ 
&&+\int_0^t\bar{\epsilon}_h(\hat{f},\theta)\;ds+\left[(u_{t}(0),\Pi_h^*\hat \theta)-(u_{h,t}(0),\Pi_h^*\hat \theta)_h\right] \nonumber \\
&=& \frac{1}{2}\|\theta(0)\|_h^2+ \int_0^tI(s)\,ds+J_1+J_2+J_3+J_4+J_5.
\label{ap8}
\end{eqnarray}
Note that estimates for the first two terms on the right hand sides of (\ref{ap8})
have already been derived in  Theorem \ref{TH1}.
For $J_1$,  use the definition of $\hat{S}$ and integrate by parts to arrive at
\begin{eqnarray} 
J_1&=& \int_{0}^{t}\beA(\hat{\theta},\theta)\;ds+
\int_{0}^{t} \int_{0}^{s}\beB(\tau,\tau;\hat{\theta}(\tau),\theta(s))\;d\tau \;ds\\
&&-\int_{0}^{t}\int_{0}^{s}\int_{0}^{\tau}
\bar\epsilon_{B_{\tau'}}(\tau,\tau';\hat{\theta}(\tau'),\theta) \;d\tau' \;d\tau \;ds\\
&=& J_{11}+J_{12}+J_{13}.
\end{eqnarray}
For $J_{11}$, a use of Lemma \ref{L1} with the inverse inequality (\ref{inv}) yields
\begin{equation}
|J_{11}|\leq   \int_0^t|\bar{\epsilon}_A(\hat{\theta},\theta)|ds \leq Ch^2\int_0^t
\|\theta\|_1\|\hat{\theta}\|_1ds \leq C h\;\int_0^t\|\theta\|_0\|\hat{\theta}\|_1\;ds.
\label{ap13}
\end{equation}
For $J_{12}$, an  integration by parts shows
\begin{eqnarray*}
|J_{12}|&=&\left|\int_0^t \bar\epsilon_B(s,s;\hat{\theta}(s),\hat{\theta}(t))ds-\int_0^t
\bar\epsilon_B(s,s;\hat{\theta}(s),\hat{\theta}(s))ds\right|\\
&\leq & Ch^2\left\{\|\hat{\theta}(t)\|_1\int_0^t\|\hat{\theta}(s)\|_1ds+
\int_0^t\|\hat{\theta}(s)\|_1^2ds\right\}.
\end{eqnarray*}
Similarly for $J_{13}$, we have
\begin{eqnarray*}
|J_{13}|&=& \left|\int_0^t\int_0^s \bar\epsilon_{B_{\tau}}(s,\tau;\hat{\theta}(\tau),\hat{\theta}(t))d\tau ds
-\int_0^t\int_0^s \bar\epsilon_{B_{\tau}}(s,\tau;\hat{\theta}(\tau),\hat{\theta}(s))d\tau ds\right|\\
&\leq & 
C(T)h^2\left\{\|\hat{\theta}(t)\|_1\int_0^t\|\hat{\theta}(s)\|_1ds+
\int_0^t\|\hat{\theta}(s)\|_1^2ds\right\}.
\end{eqnarray*}
For $J_2$, we obtain
\begin{equation}
|J_{2}|\leq  \|\bpp(\rv u)\|_{-1,h}\|\hat{\theta}\|_1+
\int_0^t\|\bpp_s(\rv u)\|_{-1,h}\|\hat{\theta}\|_1\;ds.
\label{ap10}
\end{equation}
To bound  $J_3$ and $J_4$, we integrate by parts and apply Lemma \ref{L21} to arrive at  
\begin{eqnarray}
|J_3|&\leq & |\bar{\epsilon}_h((\rv u)_t,\hat{\theta})|+\int_0^t|\bar{\epsilon}_h((\rv u)_{tt},\hat{\theta})|ds \nonumber \\
&\leq & Ch^2 \left((\|\rho_t\|_1+\|u_t\|_1)\|\hat{\theta}\|_1 + 
\int_{0}^t (\|\rho_{tt}\|_1+\|u_{tt}\|_1)\|\hat{\theta}\|_1\;ds\right)
\end{eqnarray}
and
\begin{eqnarray}
|J_{4}|&\leq& |\bar{\epsilon}_h(\hat{f},\hat{\theta})|+\int_0^t|\bar{\epsilon}_h(f,\hat{\theta})|ds \nonumber \\ &\leq&
Ch^2\Big(\|\hat{f}\|_2\|\hat{\theta}\|_1+\int_0^t\|f\|_2\|\hat{\theta}\|_1ds\Big).
\end{eqnarray}
Finally, since  $u_{h,t}(0)=\Pi_h u_t(0),$ we have 
$J_{5}=(u_{t}(0)-\Pi_h u_t(0),\Pi_h^*\hat{\theta})+\bar{\epsilon}_h(\Pi_h
u_{t}(0),\hat{\theta})$. Hence,

\begin{eqnarray}
|J_{5}|&\leq&
Ch^2\left(\|u_t(0)\|_2+\|\Pi_h u_t(0)\|_1\right)\;\|\hat{\theta}\|_1\nonumber \\
&\leq& C h^2
\|u_t(0)\|_2\;\|\hat{\theta}\|_1\leq C h^2
\|u_t(0)\|_2\|\hat{\theta}\|_1.
\label{ap12}
\end{eqnarray}
Substitute (\ref{ap13})-(\ref{ap12}) in (\ref{ap8}). We use the coercivity property of the bilinear 
form $A(\cdot,\cdot)$ and the equivalence of  the  norms, and proceed as in   
Theorem \ref{TH1} to 
complete the rest of the proof.\hfill{{\rule{2.5mm}{2.5mm}}}


Finally, we prove  quasi-optimal maximum norm estimate.
\begin{theorem}
\label{thnew}
Let $u$ and $u_h$  be the solutions of $(\ref{a})$ and $(\ref{c}),$ respectively. Further, 
let the assumptions of Lemma~$\ref{H3}$ hold. Then, 
$$
\|u(t)-u_h(t)\|_\infty \leq C(T)h^2 \left(\log\frac{1}{h}\right) \Big(
\|u_0\|_4+\|u_1\|_3+\|D_t^3 f\|_{L^1(L^2)}
+\sum_{j=0}^{2} \|D_t^j f\|_{L^1(H^{2-j})}\Big), 
$$
where $C(T)$ is a positive constant, independent of $h.$
\end{theorem}
\noindent Proof: Since $u_h(0)=\rv u_0$, it follows that $\theta(0)=0.$  Then, we  modify our 
estimates for $J_2$ to $J_4$ in (\ref{asp1}) to arrive at a superconvergence result 
for $\theta$ in $H^1$- norm 
\begin{eqnarray}
\|\theta_t\|_0+\|{\theta}\|_1&\leq& C(T)h^2 \Big(\|u_0\|_4+\|u_1\|_3+
\|D_t^3 f\|_{L^1(L^2)}\nonumber \\&&+
\sum_{j=0}^{2} \|D_t^j f\|_{L^1(H^{2-j})}\Big).
\label{NEW}
\end{eqnarray}
Now, a use  of (\ref{NEWA1}) and  (\ref{NEW}) completes the rest  of the proof.\hfill{{\rule{2.5mm}{2.5mm}}}


\section{Numerical Experiment}
\se
In this section, we present numerical results to illustrate the
performance of the finite volume element method applied to (\ref{a}).
Assume that $\cT_h$ is an admissible regular, uniform triangulation of
$\overline{\Omega}$ into closed triangles and  $0=t_0<t_1<\cdots t_M=T$ 
is a given partition of the time interval
$(0,T]$ with step length $k=\frac{T}{M}$ for some positive
integer $M$. With $U^n$ denoting the approximation of $u_h$ at $t=t_n,$
consider the discrete-in-time scheme derived in Section 5, with 
discrete $L^2$ inner product $(\cdot,\cdot)_h$ and the bilinear forms $A_h(\cdot,\cdot)$  and $B_h(t,s;\cdot,\cdot)$ 
evaluated using numerical quadrature formulae.

Thus, the  time discretization scheme is to seek $ U^n\in U_h$
for given $U^0$, such that 
\begin{eqnarray}
\frac{2}{\dt}(\ph U^{{1/2}},\Pi_h^*\chi)_h&+&\tAh(U^0,\Pi_h^*\chi)=
(f^0+\frac{2}{\dt}u_1,\Pi_h^*\chi)_h, \label{9-1}\\
(\ps U^{n},\Pi_h^*\chi)_h &+&\tAh(U^n,\Pi_h^*\chi)+
k\sum_{j=0}^{n-1}\tBh(t_n,t_{j+1/2};U^{j+1/2},\Pi_h^*\chi)\nonumber\\
&=&(f^n,\Pi_h^*\chi)_h,\label{9-1c}
\end{eqnarray}
$n\geq 1$, for all $\chi\in U_h$. The method is explicit in time in the sense that  the
calculation of $U^n$ involves only the inversion of a mass-type matrix associated with the
space $U_h$ and the corresponding dual volume element space $U_h^*$.

Let $\{\phi_j\}_{j=1,2,\cdots ,N}$ be the standard nodal basis  functions  for the 
trial space $U_h$ and  $\{\chi_j\}_{j=1,2,\cdots ,N}$ be the characteristic basis functions 
corresponding to the control volumes for the test space 
$U_h^*$. Then, express $U^n$  as 
$$\displaystyle U^n=\sum_{j=1}^N{ \alpha}_j^n \phi_j(x), \;\;\mbox{where}\;{\alpha}_j^n=U^n(x_j).$$ 
Define now the following matrices
$$
\mathbb{M}=[(\phi_i,\chi_j)_h]_{N\times N},\quad
 \mathbb{A}=[\tAh(\phi_i,\chi_j)]_{N\times N},\quad
 \mathbb{B}(t,s)=[\tBh(t,s;\phi_i,\chi_j)]_{N\times N},
$$
and the vector $\mathbb{F}(t)=[(f(t),\chi_j)_h]_{1\times N}$.
Then, for instance,  ({\ref{9-1c}}) can be written as the following  system of linear equations
 which can be solved for $\bar{\alpha}^{n+1}$: 
$$ \mathbb{M}\bar{\alpha}^{n+1}=(2\mathbb{M}-k^2\mathbb{A})\bar{\alpha}^{n-1}
-\mathbb{M}\bar{\alpha}^{n-1}-k^3\sum_{l=0}^{n-1}
\mathbb{B}(t_n,t_{l+1/2})\bar{\alpha}^{l+1/2}+k^2\mathbb{F}^n,
$$
where $\bar {\alpha}^n=\left(\alpha_1^n,\;\alpha_2^n,\cdots,\alpha_N^n\right)^T.$ Since we have used
mass lumping for $(\cdot,\cdot)_h,$ the mass matrix $ \mathbb{M}$ is a diagonal matrix.

In order to  illustrate the performance of the finite volume element method for 
solving (\ref{a}), we consider the following test problems where the computational domain
 $\Omega=(0,1)\times(0,1)$ and  the final time $T=1$.\\
{\bf Example 1:}
We choose  $u_0(x,y)=\sin(\pi x)\sin(\pi y),\; u_1(x,y)=\sin(\pi x)\sin(\pi y),\; 
A=I$ and  $B(t,s)=e^{(t-s)}I$. 
The function $f$ is chosen  so that the exact solution is  
$$u=e^t\sin(\pi x)\sin(\pi y).$$

\noindent
{\bf Example 2:}
Set $u_0(x,y)=xy(x-1)(y-1),\; u_1(x,y)=xy(x-1)(y-1),\; A=\begin{pmatrix}
1+x^2 & 0\\
0 & 1+x^2
\end{pmatrix}$ and  $B(t,s)=e^{(t-s)}A$. 
The function $f$ is chosen in such a way that the exact solution is  
$$u=e^tx y (x-1)(y-1).$$

\begin{figure}
    \begin{center}
   \includegraphics*[width=10.5cm,height=7.5cm]{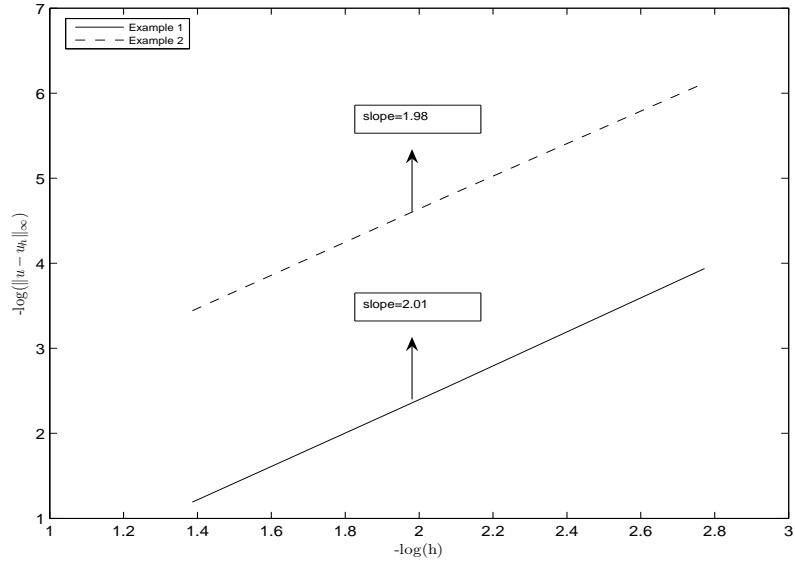}
    \caption{Convergence order estimate in $L^{\infty}$-norm.}
    \label{fig1}
    \end{center}    
\end{figure}

The order of convergence is computed in $L^{\infty}$ norm. In both examples,
Fig~\ref{fig1}  
shows that the computed order of convergence for $\|u-u_h\|_{\infty}$ in the log-log 
scale matches with the theoretical order of convergence that we have derived.

\hspace*{1cm}

{\bf { Acknowledgements.}}  
The two authors gratefully acknowledge the research support of the Department of Science and Technology, 
Government of India through the National Programme on Differential Equations: Theory, Computation and 
Applications vide DST Project No.SERB/F/1279/2011-2012, and 
the support by Sultan Qaboos University under Grant IG/SCI/DOMS/13/02.

\end{document}